\pdfoutput=1
\documentclass[12pt]{article}
\usepackage{graphicx}
\usepackage[a4paper]{geometry}
\usepackage[linktocpage=true]{hyperref}
\usepackage{float}
\usepackage{newfloat}
\usepackage{epigraph}
\usepackage{setspace}
\usepackage{xcolor}
\usepackage{tocloft}
\usepackage{graphicx}

\usepackage[font=footnotesize]{caption}
\usepackage{enumerate}
\usepackage{array}
\providecommand{\keywords}[1]{
	\par \vspace{1 em} \noindent \em \parbox{\linewidth}{\textbf{keywords} ---\ #1}}

\providecommand{\msc}[1]{
	\par \vspace{-0.55 em} \noindent \parbox{\linewidth}{\textbf{2010 MSC:} \em  #1}}

\renewcommand{\figurename}{Fig.}

\newcommand{\ipt}{
	\[y = f(x) = {x^{{x^{{x^{{x^ {\mathinner{\mkern2mu\raise1pt\hbox{.}\mkern2mu
										\raise4pt\hbox{.}\mkern2mu\raise7pt\hbox{.}\mkern1mu}} }}}}}}}\]	
}

\emergencystretch=\maxdimen
\graphicspath{{./Figures/}}
\begin{document}

\title{
\begin{spacing}{1.5}
	\vspace{-48pt}
	\LARGE{The strange properties of the infinite power tower}\\
	\vspace{-16pt}
	\textcolor{lightgray}{\hrule}
	\vspace{4pt}
	\large\textit{An ``investigative math'' approach for young students}
\end{spacing}}
\author{Luca Moroni\footnote{Liceo Scientifico "Donatelli-Pascal" - Milan - Italy}}
\date{\vspace{-15pt} \footnotesize{(August 2019)}}
\maketitle
\vspace{-32pt}
\epigraph{Nevertheless, the fact is that there is nothing as dreamy and poetic, 
nothing as radical, subversive, and psychedelic, as mathematics.} 
{Paul Lockhart -- ``\textit{A Mathematician's Lament}''}


\begin{abstract}
\setlength\parindent{0pt}
\hspace{-16pt}In this article we investigate some "unexpected" properties of the ``\textit{Infinite Power Tower}\footnote{ We adopt 
here the quite popular term \textit{power tower} even if it is not entirely correct. The 
expression should instead be described in terms of \textit{exponentiations} and not \textit{powers}.}'' function (or 
``\textit{Tetration with infinite height}''):\ipt

where the ``tower'' of exponentiations has an infinite height.\\
Apart from following an initial personal curiosity, the material collected 
here is also intended as a potential guide for teachers of 
high-school/undergraduate students interested in planning an activity of 
``\textit{investigative mathematics in the classroom}'', where the knowledge is gained through the active, creative and 
cooperative use of diversified mathematical tools (and some ingenuity).

The activity should possibly be carried on with a laboratorial style, with 
no preclusions on the paths chosen and undertaken by the students and with 
little or no information imparted from the teacher's desk.

The teacher should then act just as a guide and a facilitator. 

The infinite power tower proves to be particularly well suited to 
this kind of learning activity, as the student will have to face a 
challenging function defined through a rather uncommon infinite recursive 
process. They'll then have to find the right strategies to get around the 
trickiness of this function and achieve some concrete results, without the 
help of \textit{pre-defined} procedures. 

The mathematical requisites to follow this path are: functions, properties 
of exponentials and logarithms, sequences, limits and derivatives. The 
topics presented should then be accessible to undergraduate or ``advanced 
high school'' students.
\keywords {infinite power tower, tetration, fixed-point, recursion, recursive sequence, cobweb, Euler, Lambert, Lagrange}
\msc {97A30, 00A69}
\end{abstract}
\newpage 
\setlength{\parskip}{.0 em}
\small \setcounter{tocdepth}{1}
\begin{center}
{
	\renewcommand{\baselinestretch}{0.65}\small
	\tableofcontents
	\renewcommand{\baselinestretch}{1.0}\normalsize
}
\end{center}
\setlength{\parskip}{.6 em}
\normalsize
\vspace*{1\baselineskip}


\section{Overview}
\label{sec:overview}
After presenting the infinite power tower function, its definition and its 
unexpected properties (section \textbf{\ref{sec:introduction}} - \textit{Introduction}), we start an investigation about its mathematical 
characteristics. In section \textbf{\ref{sec:generalization}} (\textit{Generalization}) we introduce the function $y=x^y$  and its inverse function $x = y^{1/y}$  that prove to be useful to give some promising clues on the infinite power tower. In section \textbf{\ref{sec:conver}} (\textit{The problem of convergence}) we introduce the problem of the convergence of the recursive sequence   leading to the infinite power tower. This problem is furtherly investigated in the sections \textbf{\ref{sec:fixed1}} (\textit{Fixed points and convergence criteria (in general)}), \textbf{\ref{sec:fixed2}} (\textit{Fixed points and convergence of the power tower (the algebraic route})) and \textbf{\ref{sec:fixed3}} (\textit{Fixed points and convergence of the power tower (the graphical route)}), where the investigation is brought forward with both algebraic and graphical methods. Section \textbf{\ref{sec:outside}} (\textit{Outside the convergence interval}) explores what happens outside the convergence interval and the emergence of a periodic cycle for the values given by the power tower function. Lastly, in section \textbf{\ref{sec:history}} (\textit{Some history about the power tower}) we briefly discuss some very interesting historical aspects on the origin of the interest about the infinite power tower (where the main characters are Lambert, Euler and Lagrange).

\newpage 

\section{Introduction}
\label{sec:introduction}
Let's define the ``\textit{Infinite Power Tower}'' function (or ``\textit{Tetration with infinite heights}'') as:\footnote{ In the rest of this 
article we'll assume that $x,y\in \Re $ and $x>0$, $y>0$ }
\ipt
where the tower of exponentiations has an infinite height.

People aware of the \textit{explosive} nature of exponential functions will guess that, if 
$x>1$, the $f\left( x \right)$ previously defined will soon blow up to 
infinity as the height of the tower is increased. But, contrary to this 
initial guess, some trial with a pocket calculator suggests that there might 
be a stable behavior for some set of values, even with $x>1$.

In fact, some numerical experiments show that if we set $x=\sqrt 2 $ , then 
\[
\sqrt 2^{\sqrt 2^{\sqrt 2^{\sqrt 2^{\sqrt 2 
^{\,..\,\,^{\,..\,\,^{\,..\,\,}}}}}}}\to 2
\]
The reason for that can be gained by the following reasoning. Since the 
sequence of exponentials is infinite, adding (or removing) one element to 
an infinite sequence shouldn't change its overall effect (like adding or subtracting a finite number to infinite). We can then follow the passages outlined below:
\[
y=\sqrt 2^{\left\{ {\sqrt 2^{\sqrt 2^{\sqrt 2^{\sqrt 2^{\sqrt 2 
^{...}}}}}} \right\}=y}\to y=\left( {\sqrt 2 } \right)^{y}\to y^{1 
\mathord{\left/ {\vphantom {1 y}} \right. \kern-\nulldelimiterspace} 
y}=\sqrt 2 =2^{1 \mathord{\left/ {\vphantom {1 2}} \right. 
\kern-\nulldelimiterspace} 2}\to y=2
\]
We could be tempted to extend and generalize the procedure in the following 
way
\[
y=x^{\left\{ {x^{x^{x^{x^{x^{x^{...}}}}}}} \right\}=y}\to y=x^{y}\to y^{1 
\mathord{\left/ {\vphantom {1 y}} \right. \kern-\nulldelimiterspace} y}=x
\]
so that, setting $y=3$ it would be $x=3^{1 \mathord{\left/ {\vphantom {1 3}} 
\right. \kern-\nulldelimiterspace} 3}=\sqrt[3]{3}$ and setting $y=4$ it 
would be $x=\sqrt[4]{4}=\sqrt 2 $.

But here we have a problem: if we set $x=\sqrt 2 $ what will we get for the 
$y$: 2 or 4?
\[
\sqrt 2^{\sqrt 2^{\sqrt 2^{\sqrt 2^{\sqrt 2 
^{\,..\,\,^{\,..\,\,^{\,..\,\,}}}}}}}=\sqrt[4]{4}^{\sqrt[4]{4}^{\sqrt[4]{4}^{\sqrt[4]{4}^{\sqrt[4]{4}^{\,..\,\,^{\,..\,\,^{\,..\,\,}}}}}}}\to 
2?\,\,\,4?
\]
Let's check again numerically.\\ 
After having defined the following recursive function in \textit{Mathematica} or in \textit{Geogebra}

\vspace{-8pt} 
\begin{tabbing}
	Mathematica:  \= \texttt{\footnotesize{PowerTower[a\textunderscore , k\textunderscore Integer] :$=$ 
			Nest[Power[a, {\#}] {\&}, 1, k]}} \\
	Geogebra: \> \texttt{\small{Iteration(a\textasciicircum x, a, n - 1)}} (where a and n can be 
	defined as sliders)
\end{tabbing}
\vspace{-5pt} we find that a tower with height$=$1000, starting from $x=\sqrt 2 $, yields 
a result of 2 (as expected, anyway not ``4''), but if the starting point is 
$\sqrt[3]{3}$ the result is not 3 (as previously supposed), but rather a 
mysterious 2.47805. 

\begin{figure}[H]	
	\includegraphics[width=60mm,scale=1]{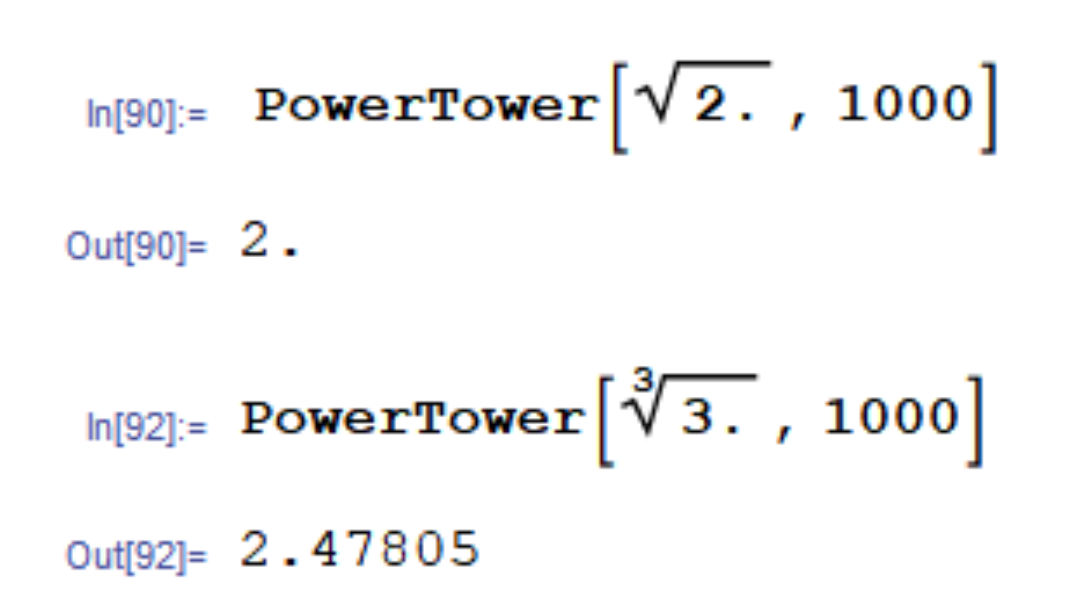}
	\label{fig000}
\end{figure}
\vspace{-18pt}
The same results are confirmed when the height of the exponentiations is 
increased to even higher values so that we can be confident that, for these 
values of the $x$, there is a definite value for the $y$.\\
Another strange thing happens when we give the $x$ some values close to 0 and 
consider odd/even numbers for the height of the tower:

\begin{figure}[H]
	\includegraphics[width=60mm,scale=1]{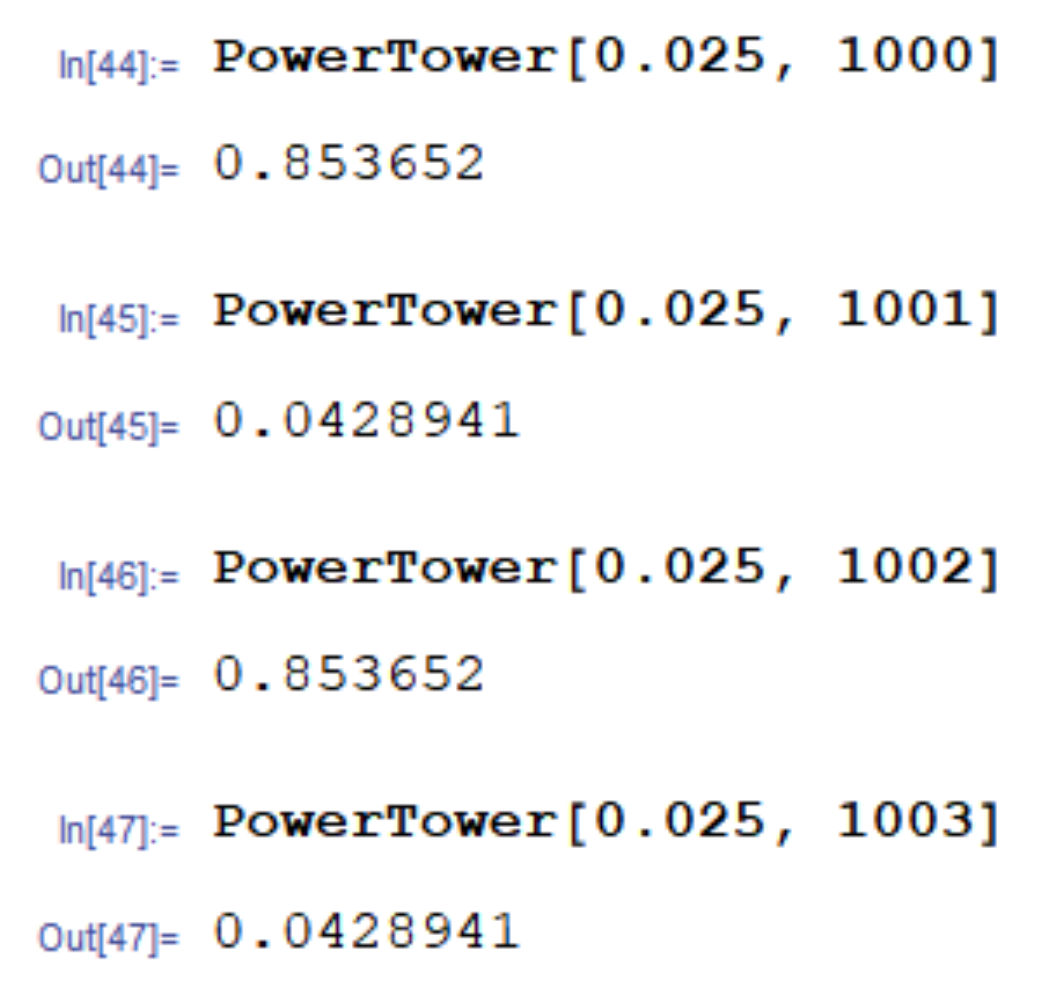}
\label{fig001}
\end{figure}
\vspace{-18pt}
It seems that a small change in the height of the tower may produce a 
relevant change in the result. How can it be?

These initial experiments suggest the following practical questions to be 
addressed: 

\begin{itemize}
\item Why is $f\left( {\sqrt 2 } \right)=2$ (and not 4)?
\item Why is not $f\left( {\sqrt[3]{3}} \right)=3$?
\item For which values of $x$ do we have a definite (finite) value of $y$?
\item Why do we sometimes get two different values whit small changes in the height of the tower?
\end{itemize}
\newpage
\section{Generalization}
\label{sec:generalization}
We have previously defined the infinite power tower function as: 
\ipt
where the tower of exponentiations has an infinite height. 

Firstly, let's make clear what is the conventional meaning of applying 
subsequent exponentiations.

In order to do so it's convenient to start from the definition of the 
related functions representing towers with \textbf{\textit{finite heights}}. 
It will be:

$y_{n} =f_{n} \left( x \right)=x^{x^{x^{x^{...}}}}$ \textbraceright $^{n\, 
times}$

so it is $y_{1} =f_{1} \left( x \right)=x,\,\,\,\,y_{2} =f_{2} \left( x 
\right)=x^{x},\,\,\,y_{3} =f_{3} \left( x \right)=x^{x^{x}}$ and so on. 

It's important to observe that it is $x^{x^{x}}=x^{\left( {x^{x}} \right)}$ 
and not $x^{x^{x}}=\left( {x^{x}} \right)^{x}=x^{x^{2}}$.

This means that 
the tower is built from the highest exponent downwards to the lowest level.

This example shows the difference between a downwards and an upwards 
construction:
\[
3^{\left(3^3   \right)}=3^{27}=7\,625\,597\,484\,987\ne \left( {3^{3}} 
\right)^{3}=9^{3}=729
\]
Using the definition of the ``finite'' tower $f_{n} \left( x \right)$, the 
infinite power tower can then be re-defined as follows:
\[
y=f(x)=x^{x^{x^{x^{{\mathinner{\mkern2mu\raise1pt\hbox{.}\mkern2mu 
\raise4pt\hbox{.}\mkern2mu\raise7pt\hbox{.}\mkern1mu}}}}}}=\lim\limits_{n\to \infty } f_{n} \left( x \right)
\]

Alternatively, we can build the sequence of functions $\left\{ {y_{1} 
,\,\,y_{2} ,\,\,y_{3} ,\,\,y_{4} ,...,\,\,y_{n} ,...} \right\}$ and take 
advantage of the fact that this sequence can be defined recursively as:
\[
\left\{ {{\begin{array}{*{20}c}
 {y_{1} =x} \\
 {y_{n+1} =x^{y_{n} }} \\
\end{array} }} \right.	f\left( x \right)=\lim\limits_{n\to \infty } y_{n} 
\]
It's easy to check that with above definition we have
\[
\begin{array}{l}
 y_{1} =x \\ 
 y_{2} =x^{y_{1} }=x^{x} \\ 
 y_{3} =x^{y_{2} }=x^{x^{x}} \\ 
 ... \\ 
 \end{array}
\]
reproducing, when $n\to \infty $, our infinite power tower.

After having clarified the meaning of the infinite power tower function 
$y=f\left( x \right)$ we can say that, if it converges to some finite value 
$y$, than it is $y=x^{\left\{ 
{x^{x^{x^{x^{{\mathinner{\mkern2mu\raise1pt\hbox{.}\mkern2mu 
\raise4pt\hbox{.}\mkern2mu\raise7pt\hbox{.}\mkern1mu}}}}}}} \right\}=y}\to 
y=x^{y}$

The inverse function will then be

\[x=g\left( y \right)=y^{1 \mathord{\left/ {\vphantom {1 y}} \right. 
\kern-\nulldelimiterspace} y}  \;\;\; (g=f^{-1})\]

Unlike $y=x^{y}$ (that's not the expression in explicit form of a function), 
this appears to be a well-defined function (although mapping $y\mapsto x)$ 
for any value $y>0$. So, let's study the characteristics of this function
$x=g\left( y \right)=y^{1 \mathord{\left/ {\vphantom {1 y}} \right. 
\kern-\nulldelimiterspace} y}$ to get some insight on the function $f\left( x 
\right)$ we are mostly interested in.

Note that we'll use the following useful identity in some calculation: 
$x=y^{1 \mathord{\left/ {\vphantom {1 y}} \right. \kern-\nulldelimiterspace} 
y}=e^{\ln y^{1 \mathord{\left/ {\vphantom {1 y}} \right. 
\kern-\nulldelimiterspace} y}}=e^{l{ny} \mathord{\left/ {\vphantom {{ny} y}} 
\right. \kern-\nulldelimiterspace} y}$

\textit{Expression}: $x=y^{1 \mathord{\left/ {\vphantom {1 y}} \right. 
\kern-\nulldelimiterspace} y}$ 

\textit{Domain}: $y>0$ 

\textit{Limits}: $\lim\limits_{y\to 0} y^{1 \mathord{\left/ {\vphantom {1 y}} \right. 
\kern-\nulldelimiterspace} y}=0^{+} \quad \lim\limits_{y\to \infty } y^{1 
\mathord{\left/ {\vphantom {1 y}} \right. \kern-\nulldelimiterspace} y}=1$

\hspace*{30pt} in fact, using the L'Hôpital's rule (H),\\
\hspace*{30pt} $y^{1/y}=e^{lny/y}$ and $\mathop {\lim }\limits_{y \to {0^ + }} \frac{{\ln y}}{y} =  - \infty $, $\mathop {\lim }\limits_{y \to  + \infty } \frac{{\ln y}}{y}\mathop  = \limits^H0$

\textit{First derivative}: $\frac{dx}{dy}=\frac{d}{dy}\left( {e^{l{ny} \mathord{\left/ {\vphantom 
{{ny} y}} \right. \kern-\nulldelimiterspace} y}} \right)=e^{l{ny} 
\mathord{\left/ {\vphantom {{ny} y}} \right. \kern-\nulldelimiterspace} 
y}\left( {\frac{1-\ln y}{y^{2}}} \right)=y^{\frac{1}{y}}\left( {\frac{1-\ln 
y}{y^{2}}} \right)$ 

\textit{Asymptotes}: the line $x=1$ is a horizontal asymptote 

\textit{Stationary points}: $\frac{dx}{dy}=0\to 1-\ln y=0\to y=e\to x=e^{1 \mathord{\left/ {\vphantom 
{1 e}} \right. \kern-\nulldelimiterspace} e}$

\hspace*{30pt} $\frac{dx}{dy}>0\to 1-\ln y>0\to y<e$ The point $M\left( {y_{M} ;x_{M} } 
\right)=\left( {e;e^{1 \mathord{\left/ {\vphantom {1 e}} \right. 
\kern-\nulldelimiterspace} e}} \right)$ is a maximum.

\textit{Second derivative: }$\frac{d^{2}x}{dy^{2}}=y^{1 \mathord{\left/ {\vphantom {1 y}} \right. 
\kern-\nulldelimiterspace} y-4}\left( {1-3y+\left( {2y+\ln y-2} \right)\cdot 
\ln y} \right)$ 

\textit{Inflection points: }$1-3y+\left( {2y+\ln y-2} \right)\cdot \ln y=0\to \left\{ 
{{\begin{array}{*{20}c}
 {y_{F1} \approx 0.5819,\,\,x_{F1} \approx 0.3944} \hfill \\
 {y_{F2} \approx 4.3678,\,\,x_{F2} \approx 1.4015} \hfill \\
\end{array} }} \right.$ 

The plot of $g\left( y \right)$ is:
\begin{figure}[H]
\centering
\includegraphics[width=125 mm,scale=1]{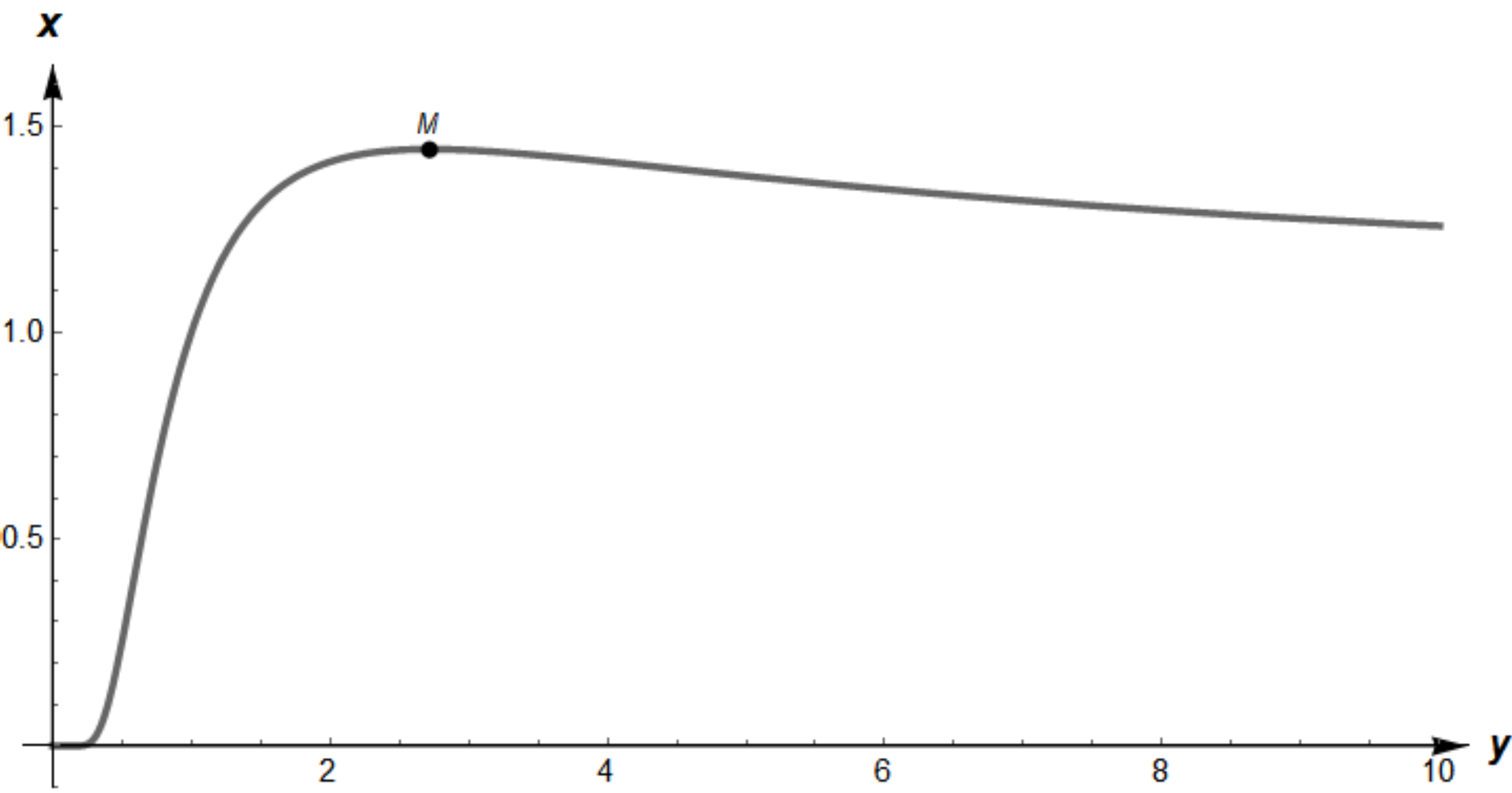}
\caption{Plot of $x=g\left( y \right)=y^{1/y}$}
\label{fig01}
\end{figure}
\vspace{-10pt}We must remember that in the plot above, differently from the usual 
conventions, the vertical axis represents the $x$ and the horizontal axis is 
the $y$. 

If we rotate the graph we get, now with the usual orientation of the axis, 
the set of points satisfying the equivalent relations:

$x=y^{1 \mathord{\left/ {\vphantom {1 y}} \right. \kern-\nulldelimiterspace} 
y}$ or $y=x^{y}$ or $y=e^{y\,\ln x}$

\begin{figure}[H]
\centering
\includegraphics[width=75 mm,scale=1]{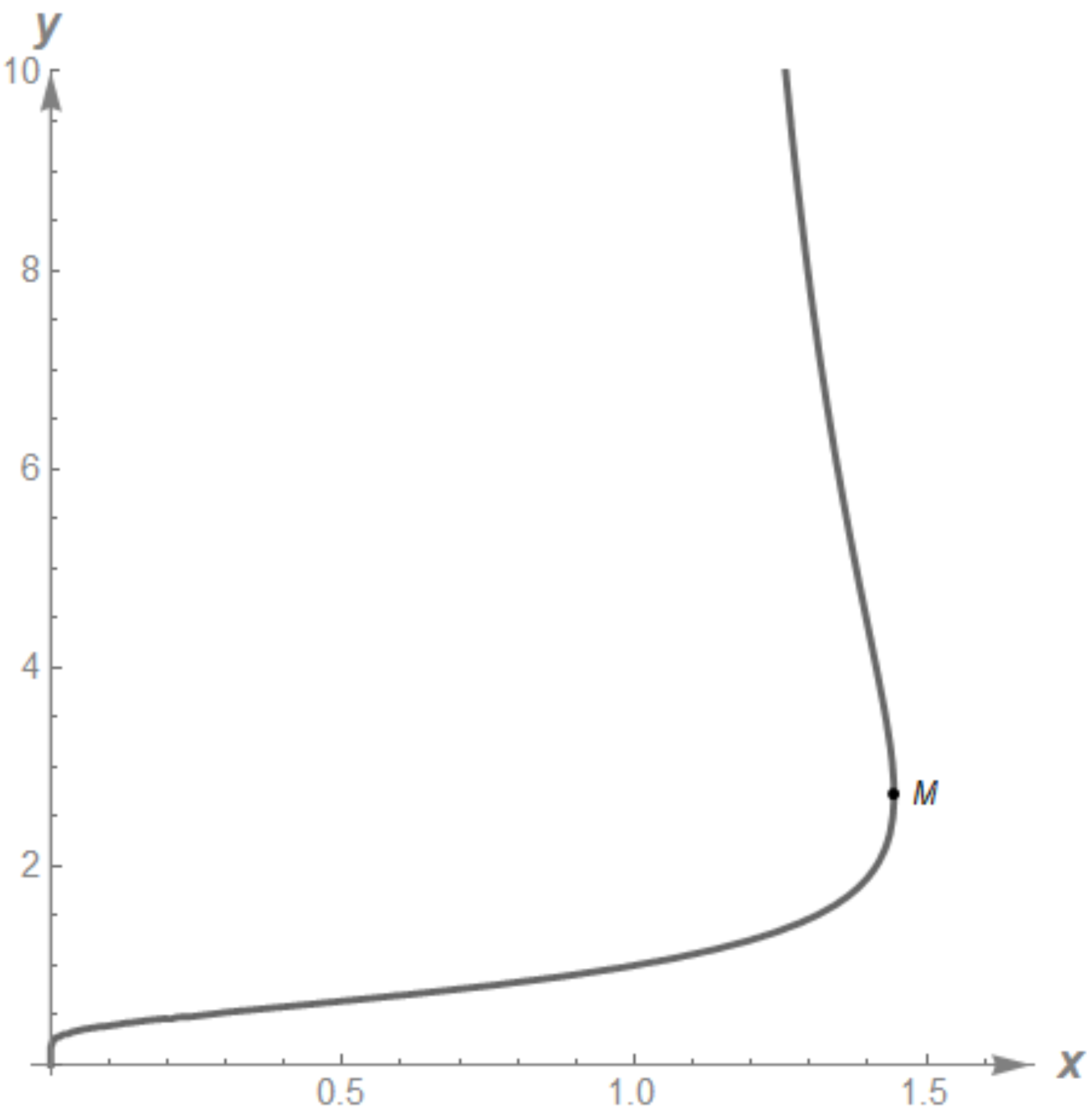}
\caption{Implicit plot of $y-x^y=0$} 
\label{fig02}
\end{figure}

But since $g\left( y \right)$ is not invertible as it is not bijective, the plot shown in \figurename~\ref{fig02} is not that of a function.

The inverse function of $g\left( y \right)$ could only be defined on a 
proper restriction of the domain of $g$, where the function is a bijection.

For example, this condition would be respected in the region defined by 
$0<x\le e^{1 \mathord{\left/ {\vphantom {1 e}} \right. 
\kern-\nulldelimiterspace} e}\wedge y>0$ and we'd get the following plot:

\begin{figure}[H]
\centering\includegraphics[width=75mm, scale=1]{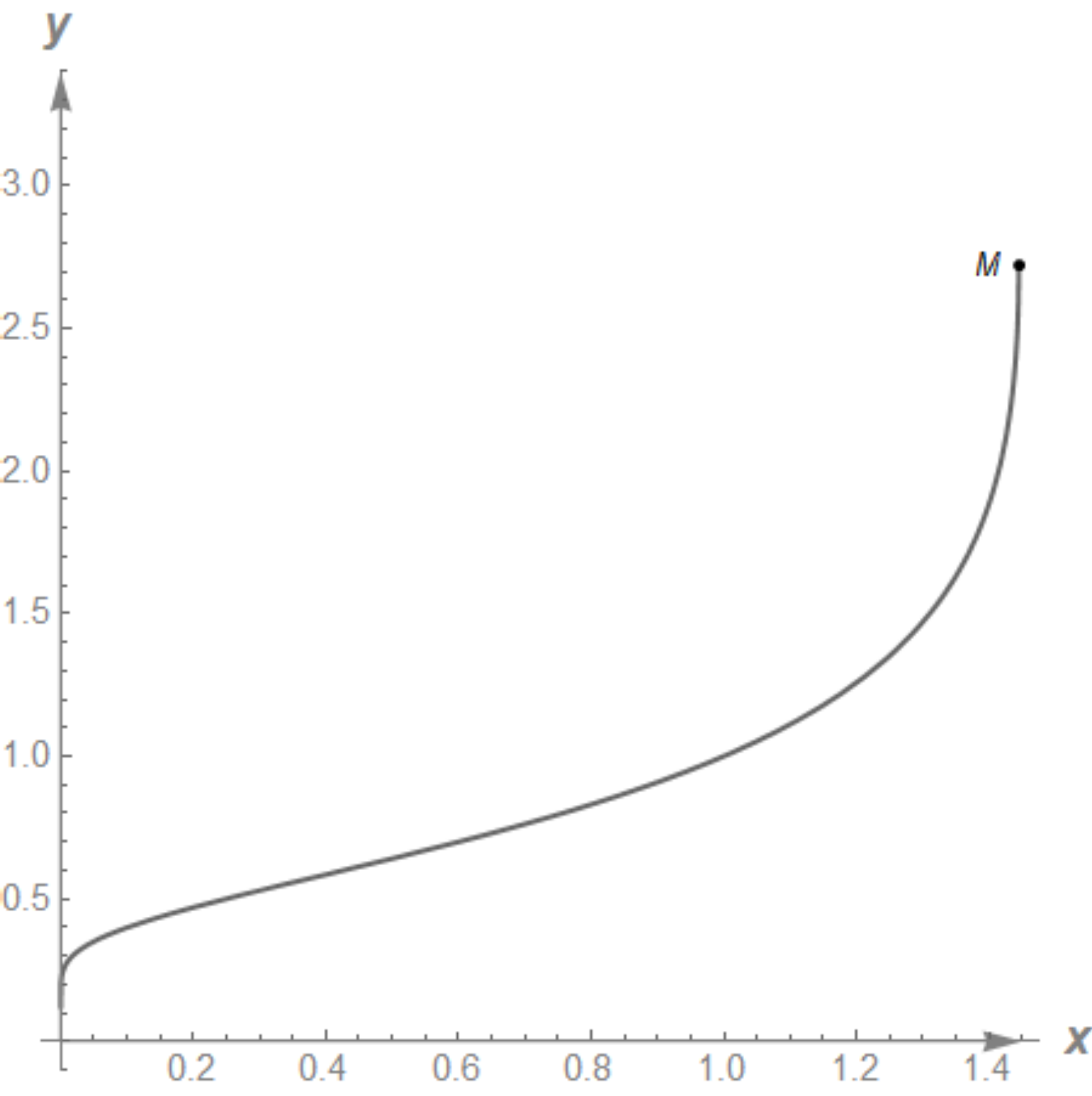}
\caption{Possible inverse function of  $g\left( y \right)$}
\label{fig03}
\end{figure}
\newpage
\section{The problem of convergence}
\label{sec:conver}
The plot of \figurename~\ref{fig02} represents all the points 
satisfying the equation $y=x^{y}$. Anyway, it would be problematic to say 
that these points are also the ones satisfying the equation of the infinite 
power tower \ipt
In fact, above equation is written in the form of a function, whilst 
$y=x^{y}$ is not the expression of a function.

Furthermore the plot tells us that for some values of the $x$ (with $1<x<e^{1 
\mathord{\left/ {\vphantom {1 e}} \right. \kern-\nulldelimiterspace} e})$ we 
would get two possible values of the $y$ and this doesn't make much sense with 
how the $f\left( x \right)$ is defined.

The problem is hidden in the following passage: 

\underline {\textbf{\textit{If the infinite power tower converges}}} than it 
is $y=x^{x^{x^{x^{x^{x^{x^{...}}}}}}}=x^{\left\{ 
{x^{x^{x^{x^{x^{x^{...}}}}}}} \right\}=y}\to y=x^{y}$ 

But the truth is that the infinite power tower \textbf{doesn't converge} for 
every values of $x$.

How can we tell that? And how can we find the interval of convergence?

We must recall that the function $f\left( x \right)$ can be defined by 
recursion as the limit of a sequence of functions with finite heights:
\[
\left\{ {{\begin{array}{ll}
 {y_{1} =x} \\
 {y_{n+1} =x^{y_{n} }} \\
\end{array} }}  \Rightarrow \right.	f\left( x \right)=\lim\limits_{n\to \infty } y_{n} 
\]

So, given some value of $x$, we can say that the sequence $\left\{ {y_{n} } 
\right\}$ converges if it stabilizes to some finite value as far as $n$ is 
increased.

In practice, the convergence requires that $\lim\limits_{n\to \infty } y_{n+1} =y_{n} $ (or $\lim\limits_{n\to \infty } 
y_{n+1} -y_{n} =0).$

To find the conditions assuring the convergence of a recursive sequence we abandon temporarily our power tower function and explore, in more general terms, sequences, fixed points and when a sequence is bound to converge to a fixed 
point. 

\section{Fixed points and convergence criteria (in general)}
\label{sec:fixed1}
In general, given a sequence defined by its starting value $y_{1} $ and by 
the recursion equation $y_{n+1} =r\left( {y_{n} } \right)$, where $r$ is a 
smooth function, a \textbf{fixed point} $y^{\ast }$ is a value satisfying 
the equation $y^{\ast }=r\left( {y^{\ast }} \right)$. The name fixed point 
means that if $y_{n} =y^{\ast }$ then $y_{n+1} =r\left( {y_{n} } 
\right)=r\left( {y^{\ast }} \right)=y^{\ast }$ and the sequence will keep on 
re-producing the same value for all future iterations.

Once we have found the fixed point(s) of a sequence by solving the equation 
$y=r\left( y \right)$ we may be interested to know if a fixed point is 
stable (or attractive) or not.

If the fixed point is attractive then, when we start close to it, we will 
end up even closer. In mathematical terms we can say that, calling $\delta 
_{n} $ the distance between $y^{\ast }$ and $y_{n} $ ($\delta_{n} >0)$ and 
starting from a point $y_{n} =y^{\ast }\pm \delta_{n} $ the subsequent term 
will be $y_{n+1} =y^{\ast }\pm \delta_{n+1} $ , and the requirement for the 
convergence is that $\delta_{n+1} <\delta_{n}$ $\forall n$. Since it is 

$\delta_{n+1} =\left| {y_{n+1} -y^{\ast }} \right|=\left| {r\left( {y_{n} } 
\right)-r\left( {y^{\ast }} \right)} \right|$ and $ \delta_{n} =\left| {y_{n} -y^{\ast }} \right| $
we have
\[
\frac{\delta_{n+1} }{\delta_{n} }=\frac{\left| {r\left( {y_{n} } 
\right)-r\left( {y^{\ast }} \right)} \right|}{\left| {y_{n} -y^{\ast }} 
\right|}
\]
The closer we are to $y^{\ast }$ the more above ratio will approximate the 
absolute value of the derivative $\left| {{r}'\left( {y^{\ast }} \right)} 
\right|$. This suggests that it's possible to use the \textit{mean value theorem} to state that there 
exist a point $\xi \in \left( {y_{n} ,y^{\ast }} \right)$ such that \[\left| 
{{r}'\left( \xi \right)} \right|=\frac{\left| {r\left( {y_{n} } 
\right)-r\left( {y^{\ast }} \right)} \right|}{\left| {y_{n} -y^{\ast }} 
\right|}\] 
In our case we can say that there exists a point $\xi \in \left( {y_{n} 
,y^{\ast }} \right)$ such that \[\frac{\delta_{n+1} }{\delta_{n} }=\left| 
{{r}'\left( \xi \right)} \right|\]

Then, if there is some interval in which it is $\left| {{r}'\left( y 
\right)} \right|<k<1 \quad \forall y\in \left( {\left| {y^{\ast }-y} 
\right|<\delta_{n} } \right)$ it will also be 
\[
\frac{\delta_{n+1} }{\delta_{n} }=\left| {{r}'\left( \xi \right)} 
\right|<k
\]
and
\[\left| {r\left( {y_{n} } \right)-r\left( {y^{\ast }} \right)} 
\right|<k\left| {y_{n} -y^{\ast }} \right|\] that is
\[\left| {y_{n+1} -y^{\ast }} \right|<k\left| {y_{n} -y^{\ast }} \right|,
\quad
\left| {y_{n+2} -y^{\ast }} \right|<k\left| {y_{n+1} -y^{\ast }} 
\right|<k^{2}\left| {y_{n} -y^{\ast }} \right|\] and so on.

We then see that if $\left| {{r}'\left( y \right)} \right|<1$ in some 
neighborhood of $y^{\ast }$ and if the starting point of the recursive 
sequence belongs to this same neighborhood, the distance to the fixed point 
reduces more and more as $n$ is increased and we'll have $\lim\limits_{n\to 
\infty } {\kern 1pt}\delta_{n} =0$ meaning that $\lim\limits_{n\to 
\infty } {\kern 1pt}y_{n} =y^{\ast }$.

In a more formal way, we state (without a complete and rigorous proof) the following 
theorem (\textit{fixed point convergence criteria}):
\begin{table}[H]
\footnotesize
\begin{center}
\begin{tabular}{|p{490pt}|}
\hline
\vspace{-4pt}
\underline {\textbf{If}} 
\begin{enumerate}[\quad 1)]
\item $r\left( y \right)$ and $r'\left( y \right)$ are continuous on $\left[ {a,b} \right]$
\item if $a\le y\le b\,\,\,\,\to \,\,\,\,a\le r\left( y \right)\le b$ (meaning that $r\left( y \right)$ is a \textit{contraction mapping})
\item $\lambda =\mathop {\max }\limits_{a\le y\le b} \left| r'\left( y \right) \right|<1$ 
\end{enumerate}
\underline {\textbf{Then}} 
\begin{enumerate}[\quad a)]
\item There exists a unique solution $y^{\ast }\in \left[ {a,b} \right]$ of the equation $y=r\left( y \right)$. 
\item For any initial starting value $y_{0} \in \left[ {a,b} \right]$ the sequence will converge to the unique fixed point: $\lim\limits_{n\to \infty } y_{n} =y^{\ast }$ 
\end{enumerate}\\
\hline
\end{tabular}
\label{tab1}
\end{center}
\end{table}
\vspace{-30pt}
The convergence/divergence character of the fixed points can be interpreted 
graphically with the so called ``\textbf{\textit{cobweb}}$''$ construction.

\begin{figure}[H]
\centering\includegraphics[width=100mm,scale=1]{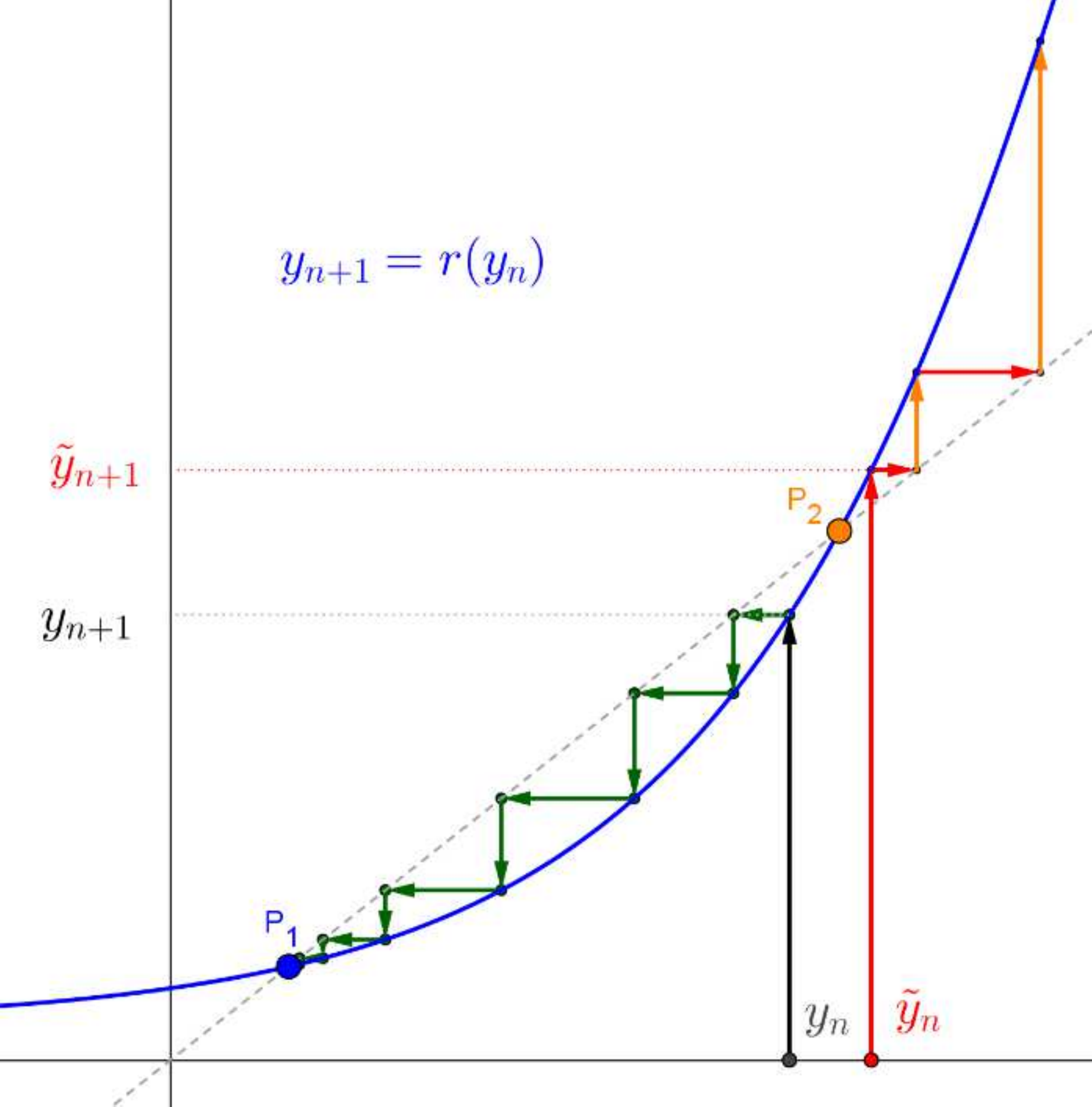}
\caption{Cobweb iteration of a sequence with an attractive ($P_1$) and a repulsive ($P_2$) fixed points.}
\label{fig04}
\end{figure}

In the following \figurename~\ref{fig04} we have the recursion 
equation $y_{n+1} =r\left( {y_{n} } \right)$ plotted with $y_{n+1} $ as a 
function of $y_{n} $. The fixed points are the intersections between 
$r\left( {y_{n} } \right)$ and the line $y_{n+1} =y_{n} $. Here we have two 
fixed points labeled $P_{1} $ and $P_{2} $. The cobweb construction shows 
that $P_{1} $ is an attractive fixed (stable) point whilst $P_{2} $ is a 
repulsive (unstable) fixed point. This is due to the fact that $\left| 
{{r}'\left( {y_{P_{1} } } \right)} \right|<1$ and $\left| {{r}'\left( 
{y_{P_{2} } } \right)} \right|>1$.

\section{Fixed points and convergence of the power tower (the algebraic route)}
\label{sec:fixed2}
In the case presented in this article we are interested in the convergence of the sequence of functions 
\[
y_{n+1} =x^{y_{n} }.
\]
Here the $x$ variable should be considered as a parameter of the recursion 
equation whose variables are the terms $y_{n} $ and $y_{n+1} $. In practice 
we have an infinite number of sequences, one for each value of $x$. 

The fixed points $\left( {y^{\ast }} \right)$ of these sequences are those 
for which it is $y_{n+1} =y_{n} $ that is those satisfying the equation
\[
y=x^{y}
\]
Using the \textit{fixed point convergence criteria} we must find the interval of the values of the $y$ (and of the $x)$ for which the first derivative of $x^{y}$ has modulus less than 1, that is
\[
\left| {\frac{d}{dy}\left( {x^{y}} \right)} \right|<1
\]
For this purpose, it's convenient to use the equivalence 

$x^{y}=e^{\ln x^{y}}=e^{y\ln x}$ and calculate the following  derivative, with respect to $y$:
\[
\frac{d}{dy}\left( {e^{y\ln x}} \right)
\]
Using the fact that for the fixed points it is $y=x^{y}$ we have
\[
\frac{d}{dy}\left( {e^{y\ln x}} \right)=e^{y\ln x}\cdot \ln x=x^{y}\cdot \ln 
x=y\cdot \ln x=\ln x^{y}=\ln y
\]
For the convergence it must then be $\left| {\ln y} \right|<1$, that is
\[
-1<\ln y<1	\to 	1 \mathord{\left/ {\vphantom {1 e}} \right. 
\kern-\nulldelimiterspace} e<y<e
\]
The corresponding values for the $x$ (since it is $x=y^{1 \mathord{\left/ 
{\vphantom {1 y}} \right. \kern-\nulldelimiterspace} y})$ are then $e^{-e}$ 
and $e^{1 \mathord{\left/ {\vphantom {1 e}} \right. 
\kern-\nulldelimiterspace} e}$ . 

Then, if we use the fixed point convergence criteria, we can say that the 
convergence is assured for 
\[
e^{-e}<x<e^{1 \mathord{\left/ {\vphantom {1 e}} \right. 
\kern-\nulldelimiterspace} e}
\]
producing stable fixed points in the interval 
\[
1 \mathord{\left/ {\vphantom {1 e}} \right. \kern-\nulldelimiterspace} 
e<y<e.
\]
\section{Fixed points and convergence of the power tower (the graphical route)}
\label{sec:fixed3}
To gain a deepest understanding of what the previous result actually means, 
we now switch to another route more rich of visual elements.

In the case of the sequence $y_{n+1} =x^{y_{n} }$ we see that the recursion 
equation is a family of exponential curves (think of the ``$x$'' as a 
parameter) and the search of the possible fixed points and their stability 
is rather simplified, mostly because these functions are strictly monotonic 
(apart the banal case with $x=1)$.

In order to simplify the notation let's rename the variables as follows:
\[
y_{n+1} =z,
\quad
y_{n} =y
\]
Then, we want to study the family of exponential functions \[z=x^{y}\] (where 
the base $x$ can be considered a parameter).

With this notation a fixed point $y^{\ast }$ is the solution of the system 
$\left\{ {{\begin{array}{ll}
 {z=x^{y}} \\
 {z=y} \\
\end{array} }} \right.$ leading to the equation $y=x^{y}$.

\newpage The character of the exponential is determined by the value of its base $x$:
\begin{itemize}
\item $x>1:$ the exponential is increasing;
\item $x=1:$ the exponential becomes the constant line $z=1$ and the original infinite power tower function becomes $y=1^{1^{1^{{\mathinner{\mkern2mu\raise1pt\hbox{.}\mkern2mu \raise4pt\hbox{.}\mkern2mu\raise7pt\hbox{.}\mkern1mu}}}}}=1$;
\item $x<1:$ the exponential is decreasing;
\end{itemize}
The positions of the curves defined by the recursive function $z=x^{y}$ with 
respect to the \textit{identity line} $z=y$ allow us to determine the possible existence of fixed points.

\begin{figure}[H]
\centering\includegraphics[width=175mm,scale=1]{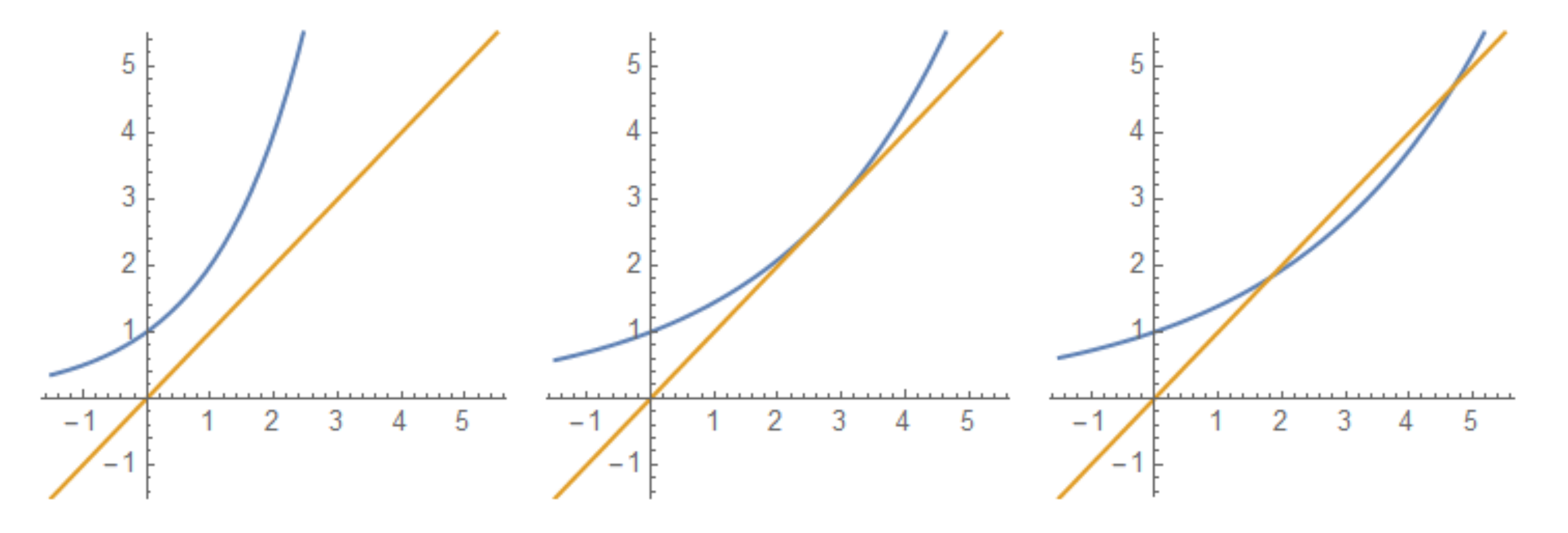}
\caption{: Possible relative positions of the exponential $z=x^{y}$ and the line 
	$z=y$ in the case $x>1$}
\label{fig05}
\end{figure}

\begin{figure}[H]
\centering\includegraphics[width=120mm,scale=1]{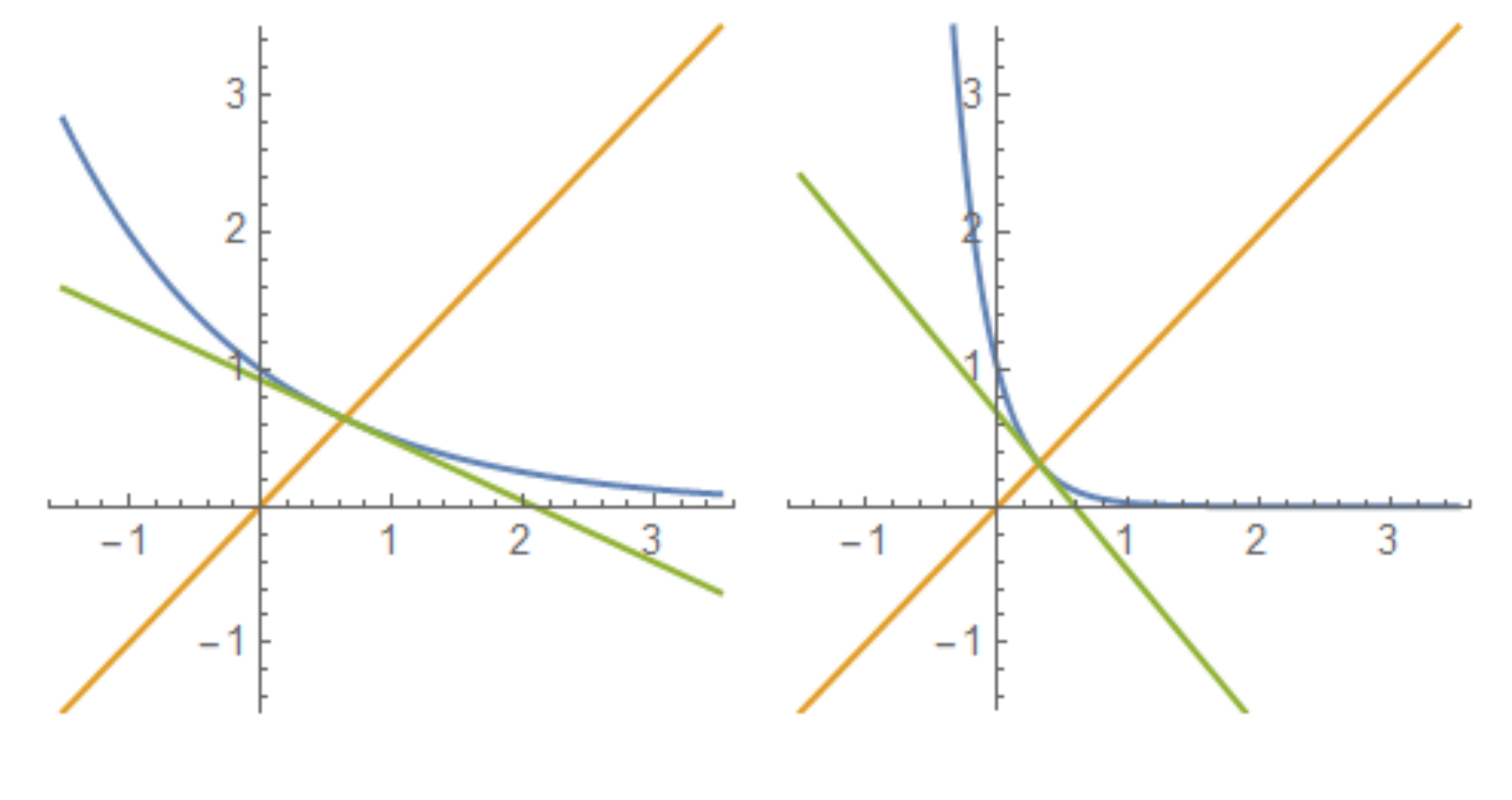}
\caption{The first derivative (slope of the tangent line) in the single fixed point in the case $0<x<1$}
\label{fig06}
\end{figure}

With $x>1$ we may have the following cases (\figurename~\ref{fig05}):

I) The exponential curve is always above the line: there are no fixed 
points.\\
II) The exponential curve is tangent to the line: there is one single fixed 
point $y_{1}^{\ast }$ (or two coincident fixed points).\\
III) The exponential curve intersects the line in two points: there are two 
distinct fixed points $y_{1}^{\ast }$ and $y_{2}^{\ast }$. 
With $x<1$ (decreasing exponential) there will always be a single 
intersection point and a single corresponding fixed point. We'll distinguish 
the following cases (\figurename~\ref{fig06}):

IV) The first derivative in the intersection point is $-1 \le \frac{d}{dy}\left( 
{x^{y}} \right)\le 0$ 

V) The first derivative in the intersection point is $\frac{d}{dy}\left( 
{x^{y}} \right)<-1$ 
Now we'll examine above 5 cases, analyze the characteristics of the fixed 
points and find which values of the ``$x$'' produce them.

\begin{figure}[H]
\centering\includegraphics[width=100mm,scale=1]{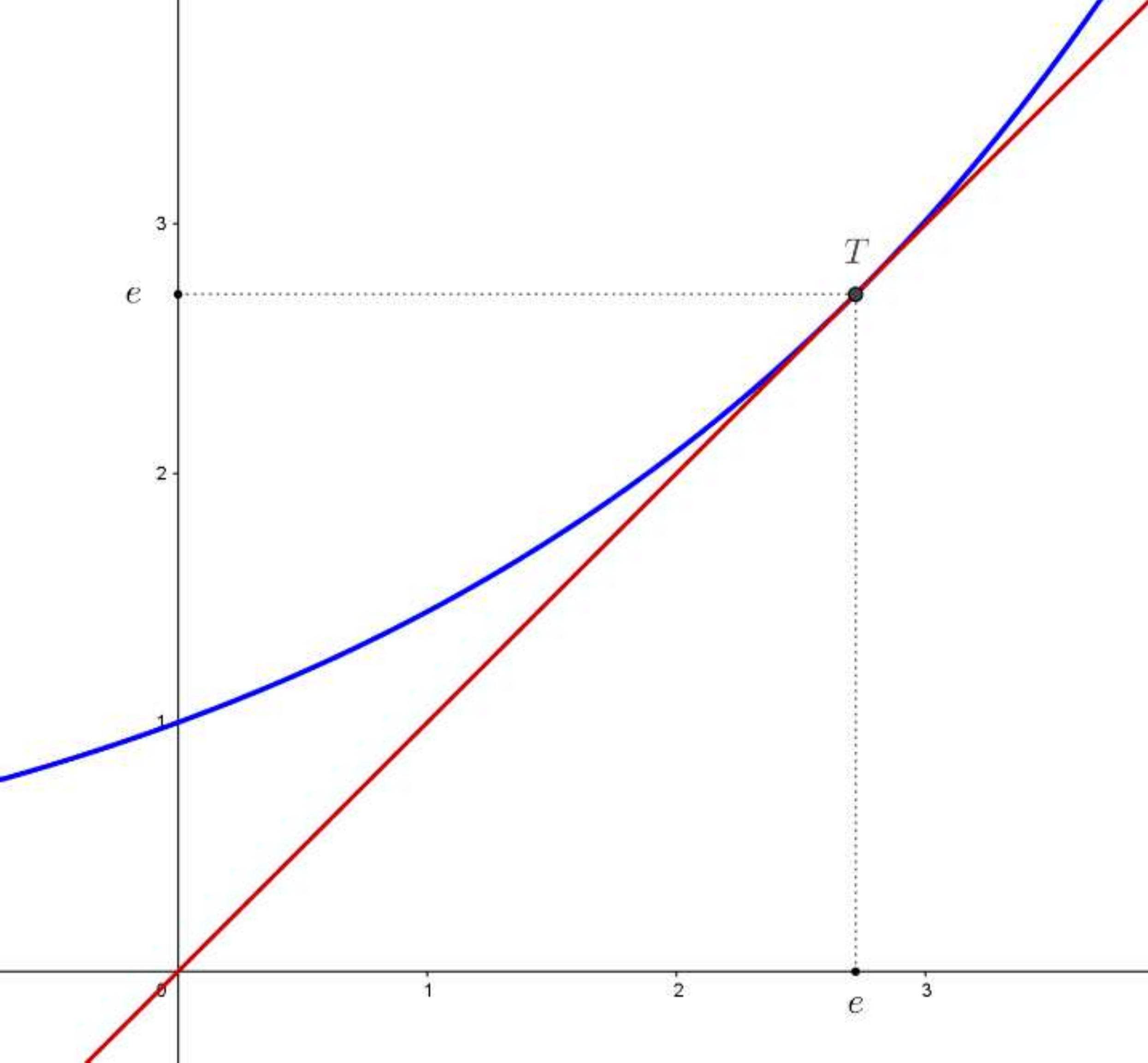}
\caption{: $x>1$ with the exponential tangent to the line $z=y$}
\label{fig07}
\end{figure}

If $x>1$ the discriminating case is that for which the exponential is 
tangent to the line (\figurename~\ref{fig07}).

So, let's look for the tangency point $T$. In this point the exponential will 
have the same slope of the line, that is $z'\left( y 
\right)=\frac{d}{dy}\left( {x^{y}} \right)=1$.\\
Then
\[
z'\left( y \right)=\frac{d}{dy}\left( {x^{y}} \right)=\frac{d}{dy}\left( 
{e^{y\ln x}} \right)=e^{y\ln x}\cdot \ln x=x^{y}\cdot \ln x=1\to y_{T} =\log 
_{x} \left( {\frac{1}{\ln x}} \right)=-\frac{\ln \ln x}{\ln x}
\]
and
\[
z_{T} =x^{y_{T} }=x^{\log_{x} \left( {1 \mathord{\left/ {\vphantom {1 {\ln 
x}}} \right. \kern-\nulldelimiterspace} {\ln x}} \right)}=\frac{1}{\ln x}
\]
We have found the point $T\left( {-\frac{\ln \ln x}{\ln x},\frac{1}{\ln x}} 
\right)$ in which it is $z'\left( y \right)=1$. But for the exponential 
curve to be tangent to the line $z=y$ we must impose that $T$ belong to that 
line, that is
\[y_{T} =z_{T} \to -\frac{\ln \ln x}{\ln x}=\frac{1}{\ln x}\to \ln \ln x=-1\to 
\ln x=e^{-1}\to x=e^{1 \mathord{\left/ {\vphantom {1 e}} \right. 
\kern-\nulldelimiterspace} e}\approx 1.445.\] 
With this value the exponential 
function becomes $z=\left( {e^{1 \mathord{\left/ {\vphantom {1 e}} \right. 
\kern-\nulldelimiterspace} e}} \right)^{y}$ and the point of tangency is 
$T\left( {e,e} \right)$.



Knowing how the base influences the graphic of a generic 
exponential curve we can also say that:

If $x>e^{1 \mathord{\left/ {\vphantom {1 e}} \right. 
\kern-\nulldelimiterspace} e}$ there's no intersection (and no fixed points 
for the recursive sequence).

If $x=e^{1 \mathord{\left/ {\vphantom {1 e}} \right. 
\kern-\nulldelimiterspace} e}$ there is a single intersection (and a single 
fixed point for the recursive sequence).

If $1<x<e^{1 \mathord{\left/ {\vphantom {1 e}} \right. 
\kern-\nulldelimiterspace} e}$ there are two intersections (and two fixed 
points for the recursive sequence).

\begin{figure}[H]
	\begin{minipage}[b]{0.45\linewidth}
		\centering
		\vspace{-40pt}
		\includegraphics[width=0.8\linewidth, scale=1]{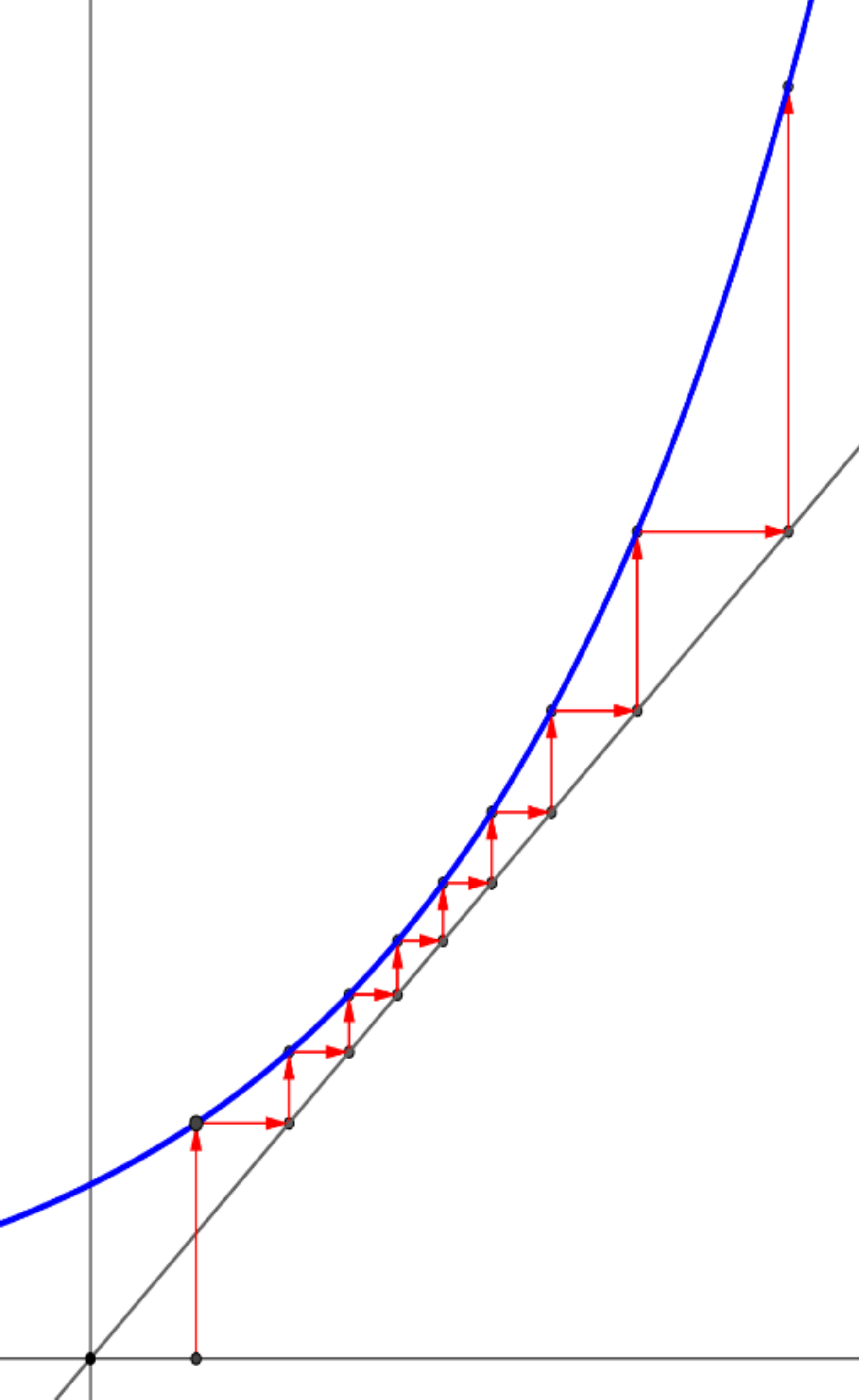}
		\caption{Cobweb diagram in the case $x>e^{1/e}$}
		\label{fig08}
	\end{minipage}\hspace*{10pt} 
	\begin{minipage}[b]{0.5\linewidth}
		\centering
		\includegraphics[width=0.8\linewidth, height=100mm]{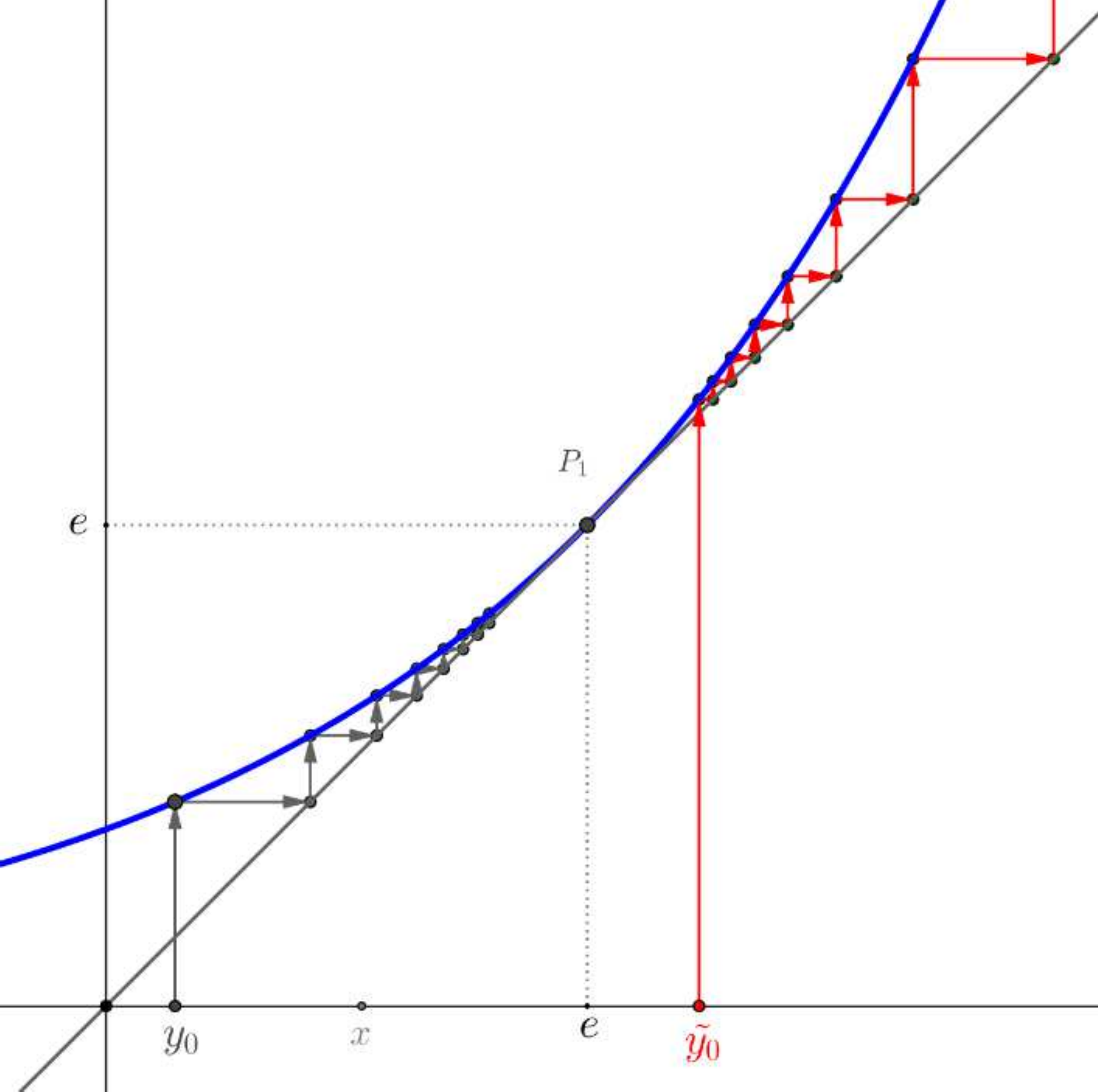}
		\caption{Cobweb diagram in the case $x=e^{1/e}$
		\label{fig09}
		}
	\end{minipage}
\end{figure}

With the cobweb diagram we can see what evolution the sequence will follow 
in these cases:

If $x>e^{1 \mathord{\left/ {\vphantom {1 e}} \right. 
\kern-\nulldelimiterspace} e}$ (\figurename~\ref{fig08}) there is 
no fixed point and, with any starting point, the sequence is bound to 
diverge to infinity.

If $x=e^{1 \mathord{\left/ {\vphantom {1 e}} \right. 
\kern-\nulldelimiterspace} e}$ (\figurename~\ref{fig09}) the cobweb 
iterations converge to $P_{1} \quad \left( {e,e} \right)$ if the starting value 
is to the left of $P_{1}$ and diverge if the starting value is to the right. 
We can say that $P_{1}$ is a ``half-stable'' saddle fixed point. Anyway, for 
the power tower sequence the starting value is $y_{0} =x=e^{1 
\mathord{\left/ {\vphantom {1 e}} \right. \kern-\nulldelimiterspace} e}$ 
that is located to the left of $P_{1} \left( {e,e} \right)$. So the sequence 
converge to $y^{\ast }=e$.

\begin{figure}[H]
\centering\includegraphics[width=100mm,scale=1]{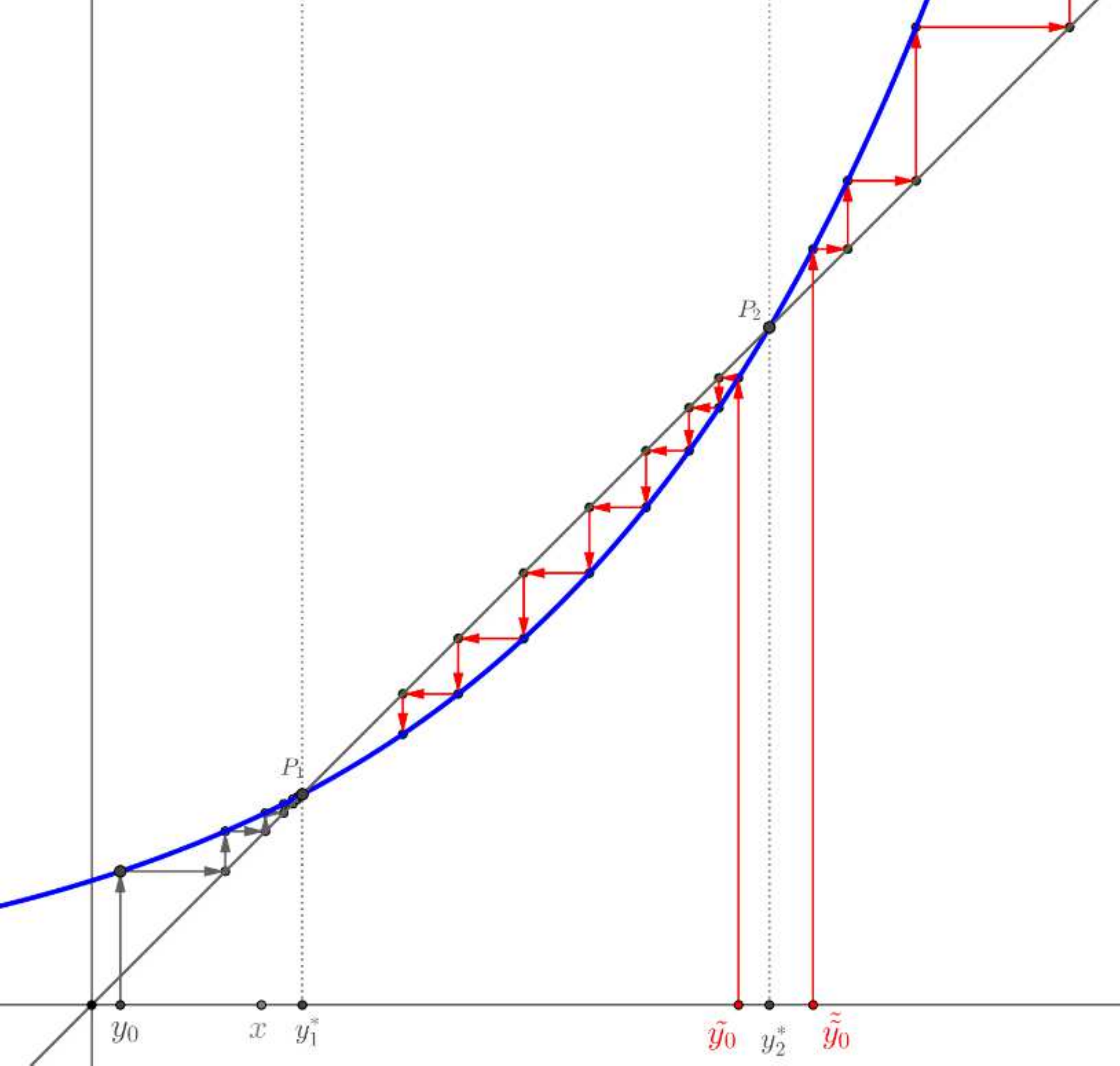}
\caption{Cobweb diagram in the case $1<x<e^{1/e}$: the fixed point $P_1$ is attractive and the fixed point $P_2$ is repulsive. The actual starting point ($x$) is in the basin of attraction of $P_1$}
\label{fig10}
\end{figure}

If $1<x<e^{1 \mathord{\left/ {\vphantom {1 e}} \right. 
\kern-\nulldelimiterspace} e}$ there are two fixed points that are the 
solutions of the of equation $y=x^{y}$. Let's call them $y^{\ast }_{1} $ and 
$y^{\ast }_{2} $(\figurename~\ref{fig10}). 

The cobweb iterations show that $y^{\ast 	 }_{1} $ is attractive and $y^{\ast 
}_{2} $ is repulsive. Furthermore the sequence will converge to $y^{\ast 
}_{1} $ for any starting point $y_{0} <y^{\ast }_{2} $ and diverge to 
infinity for $y_{0} >y^{\ast }_{2} $. Anyway, for the power tower sequence 
the starting value is $y_{0} =x$ and it's located to the left of $y^{\ast 
}_{1} $. 
In fact, for the fixed point holds the relation $y^{\ast }_{1} =x^{y^{\ast 
}_{1} }\to x=\left( {y^{\ast }_{1} } \right)^{1 \mathord{\left/ {\vphantom 
{1 {y^{\ast }_{1} }}} \right. \kern-\nulldelimiterspace} {y^{\ast }_{1} }}$ 
and if we set $x<y^{\ast }_{1} $ it must be $\left( {y^{\ast }_{1} } 
\right)^{1 \mathord{\left/ {\vphantom {1 {y^{\ast }_{1} }}} \right. 
\kern-\nulldelimiterspace} {y^{\ast }_{1} }}<y^{\ast }_{1} $. Taking the 
logarithms of both sides we have $\ln \left( {y^{\ast }_{1} } \right)^{1 
\mathord{\left/ {\vphantom {1 {y^{\ast }_{1} }}} \right. 
\kern-\nulldelimiterspace} {y^{\ast }_{1} }}<\ln y^{\ast }_{1} $ that is $1 
\mathord{\left/ {\vphantom {1 {y^{\ast }_{1} }}} \right. 
\kern-\nulldelimiterspace} {y^{\ast }_{1} }\ln \left( {y^{\ast }_{1} } 
\right)<\ln y^{\ast }_{1} \to 1 \mathord{\left/ {\vphantom {1 {y^{\ast }_{1} 
}}} \right. \kern-\nulldelimiterspace} {y^{\ast }_{1} }<1\to y^{\ast }_{1} 
>1$.

So it is $x<y^{\ast }_{1} $ if $y^{\ast }_{1} >1$. But since $z=x^{y}$ is 
increasing and it's $z\left( 0 \right)=1$, the first intersection of the 
exponential with the line $z=y$ must have a value $z>1$. This implies (since 
$y=z)$ that $y>1$. So it is $y^{\ast }_{1} >1$ and $x<y^{\ast }_{1} $. The 
sequence converges to $y^{\ast }_{1} $.

To complete our analysis let's see what happens with $0<x<1$. In this case 
the exponential curve $z=x^{y}$ is decreasing and there can be only one 
single intersection point with the line $z=y$ and a corresponding single 
fixed point. Anyway some interesting unexpected things are going to happen 
when we start analyzing the stability of that fixed point and the eventual 
convergence of the sequence to it.

\begin{figure}[H]
	\begin{minipage}[b]{0.5\linewidth}
		\centering
		\includegraphics[width=0.8\linewidth, scale=1]{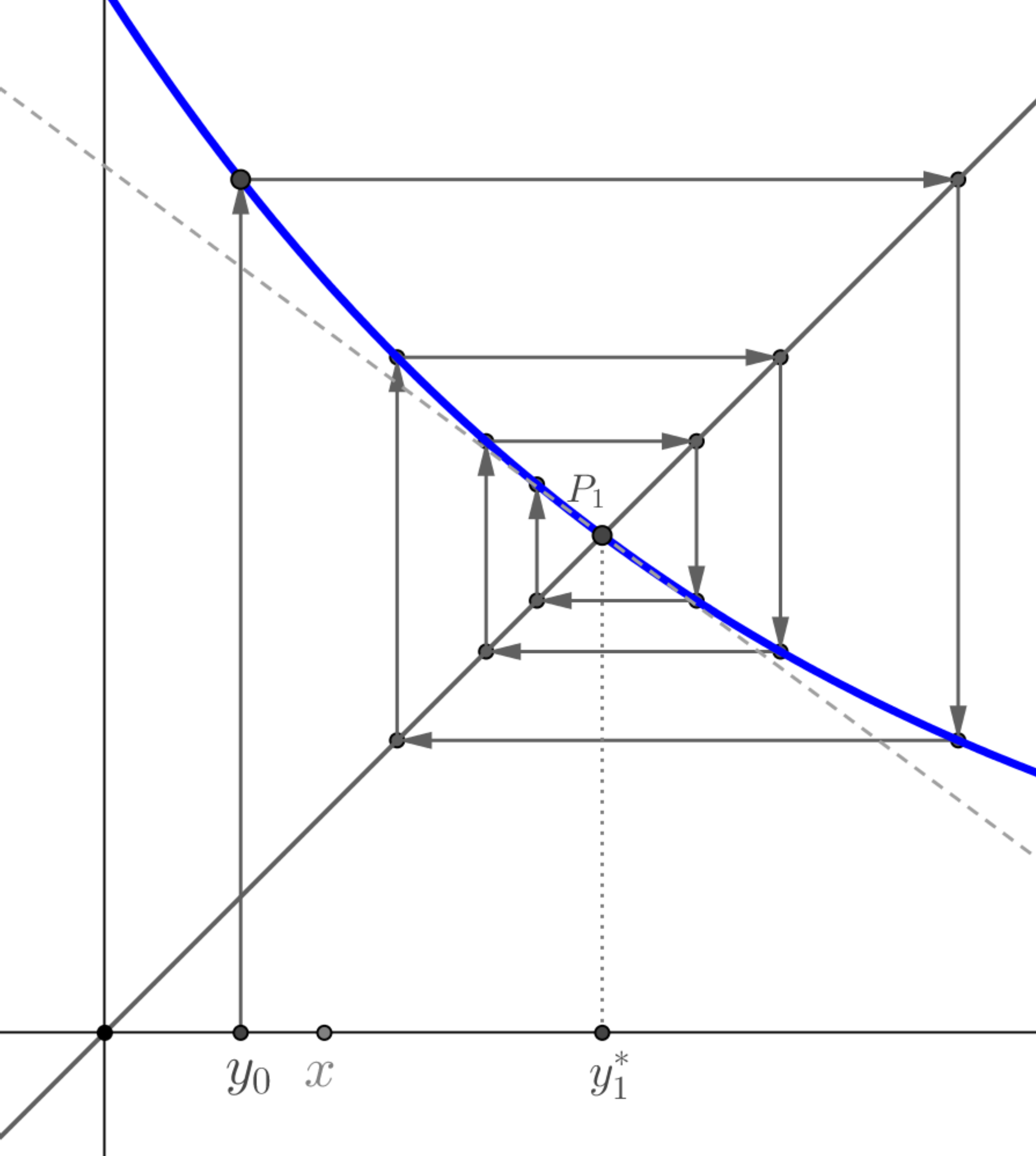}
		\caption{Cobweb iterations in the case $0<x<1$ \\with $-1 \le z'\left( {y_{0} } \right)<0$ and convergence to the fixed point.}
		\label{fig11}
	\end{minipage}\hspace*{15pt} 
	\begin{minipage}[b]{0.5\linewidth}
		\centering
		\includegraphics[width=0.8\linewidth, scale=1]{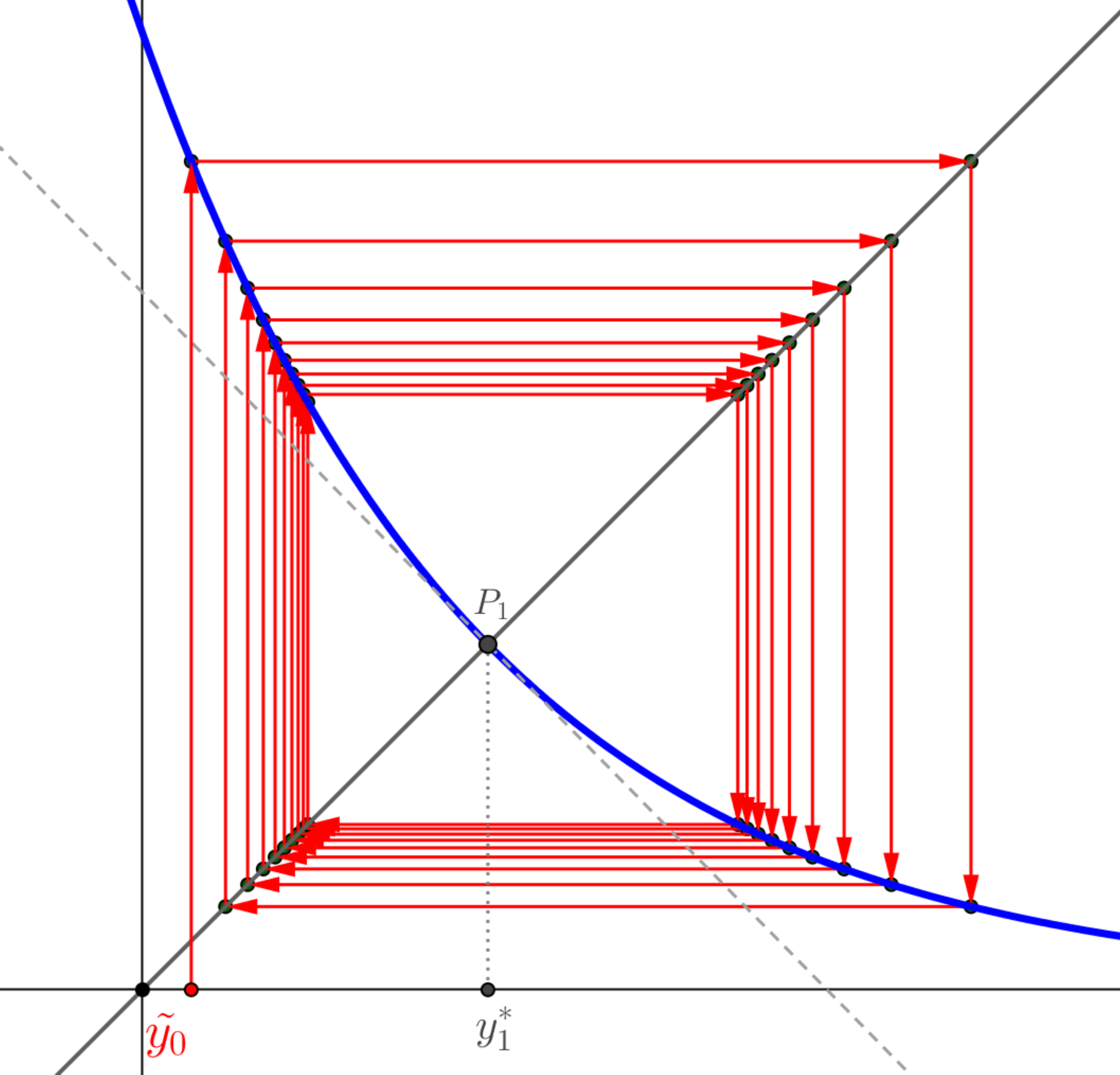}
		\caption{Cobweb iterations in the case $0<x<1$ \\with $z'\left( {y_{0} } \right)<-1$.}
		\label{fig12}
	\end{minipage}
\end{figure}

Two different cobweb iteration are presented for this case in the following 
figures, producing rather different outcomes. If the first derivative is 
$\left| {z'\left( {y^{\ast }_{1} } \right)} \right|<1$ (that is $-1 \le z'\left( 
{y^{\ast }_{1} } \right)<0$, \figurename~\ref{fig11}) the 
iterations converge to $y^{\ast }_{1} $, oscillating between values 
alternatively greater and less than that of the fixed point. We can say that 
the fixed point is attractive and that the sequence will eventually converge 
to it, whatever is the starting point.

On the contrary, if the first derivative is $\left| {z'\left( {y^{\ast }_{1} 
} \right)} \right|>1$ (that is $z'\left( {y^{\ast }_{1} } \right)<-1$, 
\figurename~\ref{fig12}) the iterations are again oscillating but 
the sequence doesn't converge to $y^{\ast }_{1} $ . Instead it stabilizes 
towards a periodic stable cycle, getting closer and closer to two alternate 
distinct fixed values.

Let's then see for what value of $x$ we have $z'\left( {y^{\ast }_{1} } 
\right)>-1$. Since we have already found that $z'\left( {y^{\ast }_{1} } 
\right)=x^{y}\ln x$ we must solve the inequality $x^{y}\ln x>-1$ with 
$y=x^{y}$ meaning $x=y^{1 \mathord{\left/ {\vphantom {1 y}} \right. 
\kern-\nulldelimiterspace} y}$. It will then be
\[
y\ln y^{1 \mathord{\left/ {\vphantom {1 y}} \right. 
\kern-\nulldelimiterspace} y}>-1\to \ln y>-1\to y>e^{-1}\to x=y^{1 
\mathord{\left/ {\vphantom {1 y}} \right. \kern-\nulldelimiterspace} 
y}>e^{-e}
\]
We can then say that the fixed point is attractive for $e^{-e}\le x<1$ and 
that we'll have a 2-cycle for $0<x<e^{-e}$.

What can we say about the two values $y_{1} $ and $y_{2} $ involved in the 
2-cycle? 

Since $y_{2} $ is the next value in the sequence after $y_{1} $ and $y_{1} $ 
is the next value in the sequence after $y_{2} $ we have $y_{2} =x^{y_{1} }$ 
and $y_{1} =x^{y_{2} }$. Let's take the power $y_{1} $ of both sides of the 
second equation to get $y_{1}^{y_{1} }=x^{y_{2} y_{1} }$. 

Inserting $x^{y_{1} }=y_{2}$ we have $y_{1}^{y_{1} }=y_{2}^{y_{2} }$ and 
$y_{1} \ln y_{1} =y_{2} \ln y_{2} $. 

Let's now say that $y_{2} $ is $p$ times $y_{1} $, that is $y_{2} =py_{1}$ and 
solve for $y_{1} $. 
\[
py_{1} \ln py_{1} =y_{1} \ln y_{1} 
\]
\[
p\left( {\ln p+\ln y_{1} } \right)=\ln y_{1} \,\,\,\to \,\,\,\left( {p-1} 
\right)\ln y_{1} =-p\ln \left( p \right)\,\,\,\to \,\,\,\ln y_{1} 
=p\frac{\ln p}{1-p}\,\,\,\to \,\,\,\
ln y_{1} =\ln \left( {p^{\frac{p}{1-p}}} 
\right)
\]
We finally have

\begin{center}
	$y_{1} =p^{\frac{p}{1-p}}$  \quad and \quad $y_{2} =py_{1} 
=p^{\frac{p}{1-p}+1}=p^{\frac{1}{1-p}}.$
\end{center}

For instance, if we set $p=2$ we have $y_{1} =1 \mathord{\left/ {\vphantom 
{1 4}} \right. \kern-\nulldelimiterspace} 4$ and $y_{2} =1 \mathord{\left/ 
{\vphantom {1 2}} \right. \kern-\nulldelimiterspace} 2$. These are the two 
values of the cycle that we'd get with $x=y_{2}^{1 \mathord{\left/ 
{\vphantom {1 y}} \right. \kern-\nulldelimiterspace} y_{1} }=y_{1}^{1 
\mathord{\left/ {\vphantom {1 y}} \right. \kern-\nulldelimiterspace} y_{2} 
}=\left( {1 \mathord{\left/ {\vphantom {1 2}} \right. 
\kern-\nulldelimiterspace} 2} \right)^{4}=\left( {1 \mathord{\left/ 
{\vphantom {1 4}} \right. \kern-\nulldelimiterspace} 4} \right)^{2}=1 
\mathord{\left/ {\vphantom {1 {16=0.0625}}} \right. 
\kern-\nulldelimiterspace} {16=0.0625}$. We can also find the two values for 
the \textit{extreme} cycle in which $y_{1} $ and $y_{2} $ have the maximum separation. 
Setting $p\to \infty $ we have
\[
y_{1} =\lim\limits_{p\to \infty } p^{\frac{p}{1-p}}\to \infty 
^{-1}=0^{+}
\]
\[
y_{2} =\lim\limits_{p\to \infty } p^{\frac{1}{1-p}}=\lim\limits_{p\to 
\infty } e^{\frac{1}{1-p}\ln p}\to e^{0^{-}}=1^{-}
\]
and, since $y_{1} =x^{y_{2} }$, $0=x^{1},\,\,\,x\to 0$.

The following table summarizes the results gained so far.

\begin{table}[H]
\footnotesize
\begin{center}
\centering
\begin{tabular}{|m{70pt}|m{100pt}|m{125pt}|m{125pt}|}
\hline
\vspace{5pt}
\textbf{Values of }\textbf{\textit{x}} & 
\vspace{5pt}
\textbf{Fixed point values} &
\vspace{5pt}
\textbf{Fixed point(s)} &
\vspace{5pt}
\textbf{Asymptotic behavior} \\[10pt]
\hline
\hline
\vspace{10pt}$x>e^{1 \mathord{\left/ {\vphantom {1 e}} \right. \kern-\nulldelimiterspace} e}$& 
	\vspace{10pt} \quad //& 
	\vspace{10pt}No fixed points& 
	\vspace{10pt}Divergence to $+\infty $ \\[14pt]
	 \hline
		 \vspace{10pt}$x=e^{1 \mathord{\left/ {\vphantom {1 e}} \right. \kern-\nulldelimiterspace} e}$& 
	\vspace{10pt}$y=e$& 
	\vspace{10pt}1 fixed point& 
	\vspace{10pt}Convergence to the f.p. \\[14pt]
	\hline
	\vspace{10pt}$1<x<e^{1 \mathord{\left/ {\vphantom {1 e}} \right. \kern-\nulldelimiterspace} e}$& 
	\vspace{10pt}$1<y<e$& 
	\vspace{10pt}2 fixed points& 
	\vspace{10pt}Convergence to the first f.p. \\[14pt]
	\hline
	\vspace{10pt}$x=1$& 
	\vspace{10pt}$y=1$& 
	\vspace{10pt}1 fixed point& 
	\vspace{10pt}Instantaneous convergence to the f.p. \\[14pt]
	\hline
	\vspace{10pt}$e^{-e}<x<1$& 
	\vspace{10pt}$1 \mathord{\left/ {\vphantom {1 e}} \right. \kern-\nulldelimiterspace} e<y<e$& 
	\vspace{10pt}1 fixed point& 
	\vspace{10pt}Convergence to the f.p. (oscillating) \\[14pt]
	\hline
	\vspace{10pt}$x=e^{-e}$& 
	\vspace{10pt}$y=1 \mathord{\left/ {\vphantom {1 e}} \right. \kern-\nulldelimiterspace} e$& 
	\vspace{10pt}1 fixed points& 
	\vspace{10pt}Convergence to the f.p. (oscillating) \\
	\hline
	\vspace{8pt}$0<x<e^{-e}$& 
	\vspace{8pt}$0<y_{1} <1 \mathord{\left/ {\vphantom {1 e}} \right. \kern-\nulldelimiterspace} e<y_{2} <1$& 
	\vspace{8pt}1 fixed point (unstable)$+$ \par a stable 2-cycle with values $y_{1} $ and $y_{2} $ & 
	\vspace{8pt}Convergence to the 2-cycle \\[14pt]
	\hline
	\vspace{8pt}$x\to 0^{+}$& 
	\vspace{0pt}$\begin{array}{l}
	y_{1} \to 0 \\ 
	y_{2} \to 1 \\ 
	\end{array}$& 
	\vspace{8pt}1 fixed point (unstable)$+$ \par a stable 2-cycle & 
	\vspace{10pt}Convergence to the 2-cycle \\[17pt]
	\hline
\multicolumn{4}{|p{215pt}|}{\vspace{1pt} Decimal values:${\begin{array}{*{20}c}
			{e\approx 2.71828} & {e^{1 \mathord{\left/ {\vphantom {1 e}} \right. \kern-\nulldelimiterspace} e}\approx 1.44467} & {1/e\approx 0.367879} & {e^{-e}=1 \mathord{\left/ {\vphantom {1 {e^{e}\approx 0.065988}}} \right. \kern-\nulldelimiterspace} {e^{e}\approx 0.065988}} \\
			\end{array} }$ } \\[17pt]
	\hline
\end{tabular}
\label{ctab01}
\caption{Convergence intervals of the infinite power tower.}
\end{center}
\end{table}
\vspace{-15pt}It's interesting to note that how the number $e$ appears in above table in 
many possible power variations.

In conclusion, we can now say that the infinite power tower converges to the 
function defined by the expression $y=x^{y}$ (or $x=y^{1 \mathord{\left/ 
{\vphantom {1 y}} \right. \kern-\nulldelimiterspace} y}$ ) for $e^{-e}\le 
x\le e^{1 \mathord{\left/ {\vphantom {1 e}} \right. 
\kern-\nulldelimiterspace} e}$ assuming values $1 \mathord{\left/ {\vphantom 
{1 e}} \right. \kern-\nulldelimiterspace} e\le y\le e$.

\begin{figure}[H]
\centering\includegraphics[width=300pt,scale=1]{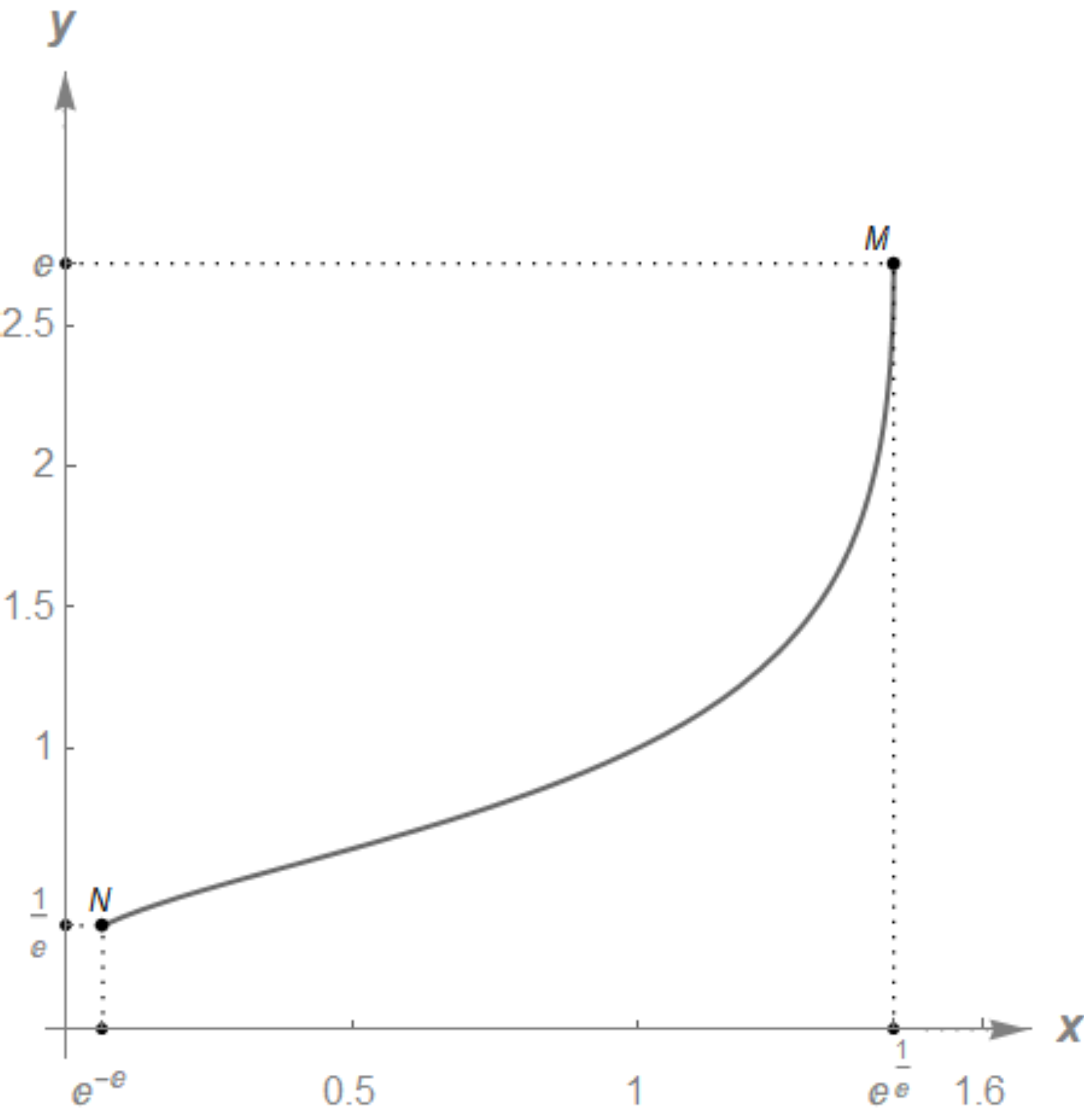}
\caption{Plot of 	$y = f(x) = {x^{{x^{{x^{{x^{\mathinner{\mkern2mu\raise1pt\hbox{.}\mkern2mu							\raise4pt\hbox{.}\mkern2mu\raise7pt\hbox{.}\mkern1mu}} }}}}}}}$	}
\label{fig13}
\end{figure}

Taking into account the information collected we can show, in 
\figurename~\ref{fig13}, the final plot of the infinite power 
tower function.

\section{Outside the convergence interval}
\label{sec:outside}
We have seen that the infinite power tower converges for $e^{-e}\le x\le 
e^{1 \mathord{\left/ {\vphantom {1 e}} \right. \kern-\nulldelimiterspace} 
e}$, assuming values $1 \mathord{\left/ {\vphantom {1 e}} \right. 
\kern-\nulldelimiterspace} e\le y\le e$.

But what happens outside the convergence interval? 

Let's try some numerical experiment with some power towers with finite (but 
rather high) height.

\begin{figure}[H]
	\centering\includegraphics[width=360pt,scale=1]{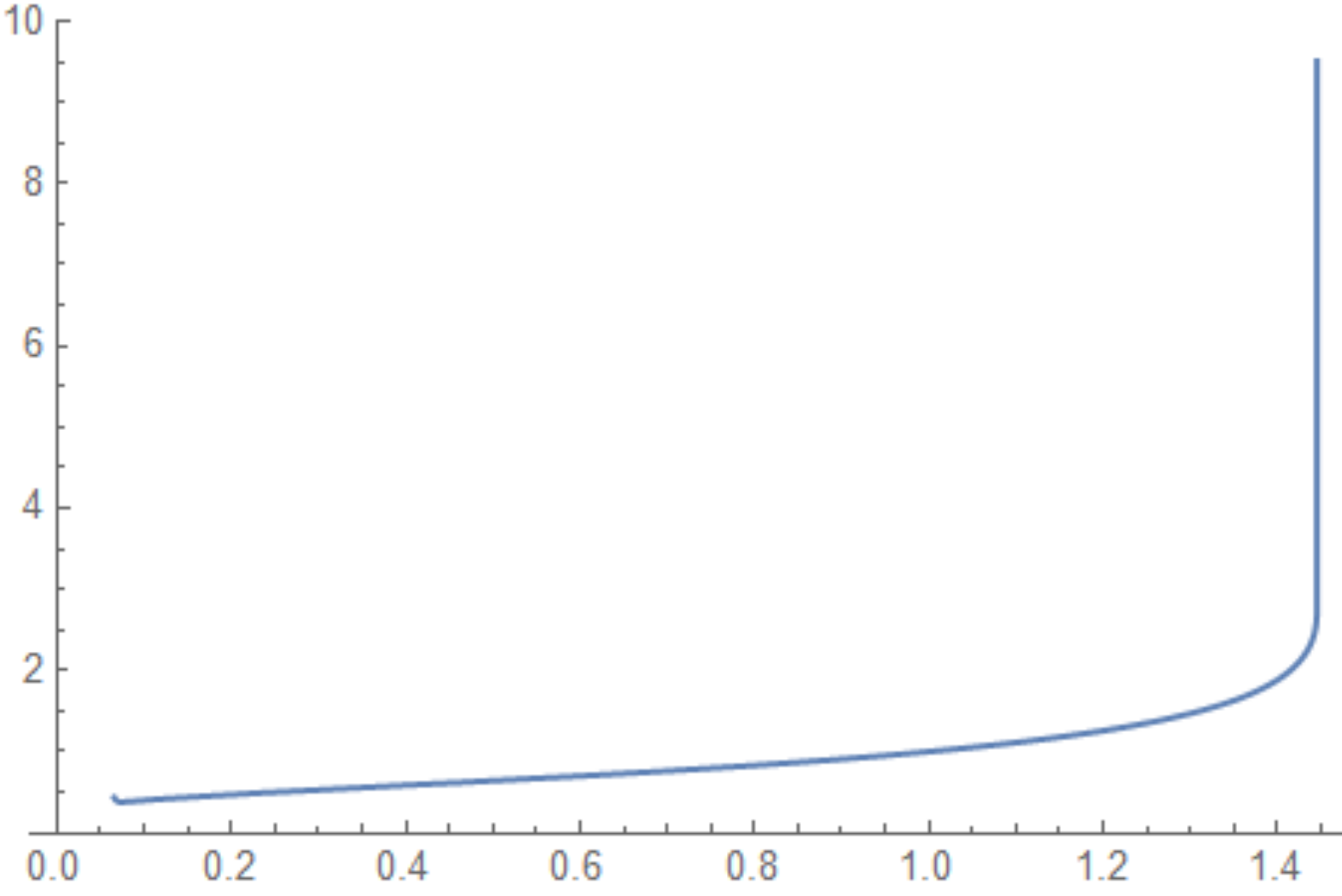}
	\caption{Plot of $y_{200} \left( x \right)=PowerTower[x,200]$}
	\label{fig14}
\end{figure}

For $x>e^{1 \mathord{\left/ {\vphantom {1 e}} \right. 
\kern-\nulldelimiterspace} e}$ the function $f\left( x \right)$ blows out 
rapidly to $+\infty $ (\figurename~\ref{fig14}). In fact we already know that there aren't fixed 
points for the sequence $y_{n+1} =x^{y_{n} }$ when we use $x>e^{1 
\mathord{\left/ {\vphantom {1 e}} \right. \kern-\nulldelimiterspace} e}$ .

\begin{figure}[H]
\centering\includegraphics[width=360pt,scale=1]{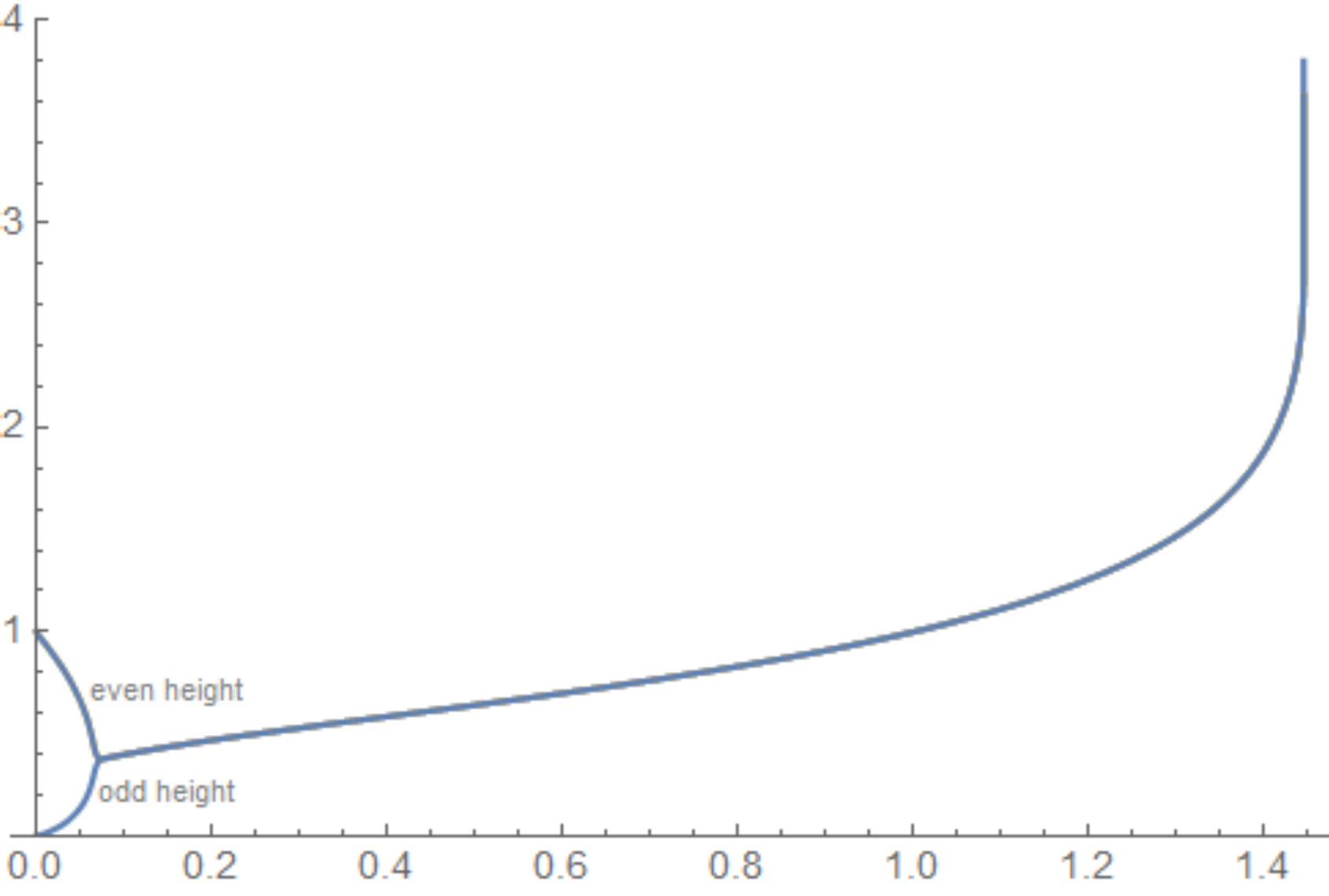}
\caption{Plot of $y_{200} \left( x \right)=PowerTower[x,200]$, $y_{201} \left( x \right)=PowerTower[x,201]$}
\label{fig15}
\end{figure}

As we have already seen, for $0<x<e^{-e}$ the sequence start to oscillate 
between two bounded values, and some numerical simulation confirms that 
behavior (\figurename~\ref{fig15}). 
The upper/lower branches of the plot correspond to an even/odd value for the 
height of the tower.

We have already seen this oscillating behavior when exploring, through the 
cobweb diagrams, the recursive sequence $y_{n+1} =x^{y_{n} }$ with 
$x<e^{-e}$. 

Let's analyze further the origin of this feature.

First we can try to calculate the limit of the finite power tower when $x\to 
0$. Let's start with $f_{2} $:
\[
\lim\limits_{x\to 0} f_{2} \left( x \right)=\lim\limits_{x\to 0} 
x^{x}=\lim\limits_{x\to 0} e^{x\ln x}
\]
we can now calculate the limit of the exponent (using the 
L'\textbf{H}\^{o}pital's rule)
\[
\lim_{ x \to 0 } x lnx=
\lim_{ x \to 0 } \frac {lnx} {1/x}\mathop  = \limits^H 
\lim_{ x \to 0 } \frac {1/x} {-1/x^2}=
\lim_{ x \to 0 }\left( -x  \right) =0^{-}
\]
so it is
\[
\lim\limits_{x\to 0} e^{x\ln x}=1^{-}
\]
Then we are on the upper branch. What changes when we increase the tower 
height?

$$\lim\limits_{x\to 0} f_{3} \left( x \right)=\lim\limits_{x\to 0} 
x^{x^{x}}=\lim\limits_{x\to 0} x^{\left( {x^{x}} \right)}$$ and since 
$\lim\limits_{x\to 0} x^{x}=1$, as seen before, the last limit has the 
form $0^{1}\to 0$

So it is $\lim\limits_{x\to 0} f_{2} \left( x \right)=\lim\limits_{x\to 0} x^{x}=1$ and $\lim\limits_{x\to 0} f_{3} \left( x 
\right)=\lim\limits_{x\to 0} x^{x^{x}}=0$.

We have a strong suspect (supported by the previous reasoning based on the 
cobweb diagrams) that these results may extend to towers with any height, 
with different values for $n$ even and $n$ odd, that is
$\lim\limits_{x\to 0} f_{2n} \left( x \right)=1$ and $\lim\limits_{x\to 
0} f_{2n+1} \left( x \right)=0$ but we can't prove this conjecture with simple tools and leave this problem to a later time.

Having observed the oscillating character of the finite power tower sequence 
for $0<x<e^{-e}$, we ask ourselves if it's possible to find the equations of 
the two distinct branches.

Calling $a$ and $b$ the two values corresponding to some $\bar{{x}}$ it must be:
\[{y_n}\left( {\overline x } \right) \to a,{\kern 1pt} {\kern 1pt} {\kern 1pt} {\kern 1pt} {y_{n + 1}}\left( {\overline x } \right) \to b,{\kern 1pt} {\kern 1pt} {\kern 1pt} {\kern 1pt} {y_{n + 2}}\left( {\overline x } \right) \to a,\,\,\,\,{y_{n + 3}}\left( {\overline x } \right) \to b, \ldots \]

This means that \textbf{the sequence built with a double recursion} should 
converge to its fixed points $a$ and $b$.

Let's see what is the form of this double recursion:
\[
\left\{ {{\begin{array}{*{20}l}
 {y_{n+1} =x^{y_{n} }} \\
 {y_{n+2} =x^{y_{n+1} }} \\
\end{array} }} \right.\Rightarrow y_{n+2} =x^{x^{y_{n} }}
\]
This sequence has stable fixed points if the derivative with respect to $y$ of 
the right side has modulus less than 1, that is
\[
\left| {\frac{d}{dy}\left( {x^{x^{y}}} \right)} \right|<1
\]
Since it is $x^{x^{y}}=e^{\ln x^{x^{y}}}=e^{x^{y}\ln x}$ the derivative to 
calculate becomes:
\[
\frac{d}{dy}\left( {e^{x^{y}\ln x}} \right)
\]
after some passage we arrive at
\[
\frac{d}{dy}\left( {e^{x^{y}\ln x}} \right)=x^{x^{y}+y}\ln^{2}x
\]
and it must be \[\left| {x^{x^{y}+y}\ln^{2}x} \right|<1\] 

Differently from before we can't find an explicit algebraic form for the 
boundary of the region of convergence.

\begin{figure}[H]
\centering\includegraphics[width=320pt,scale=1]{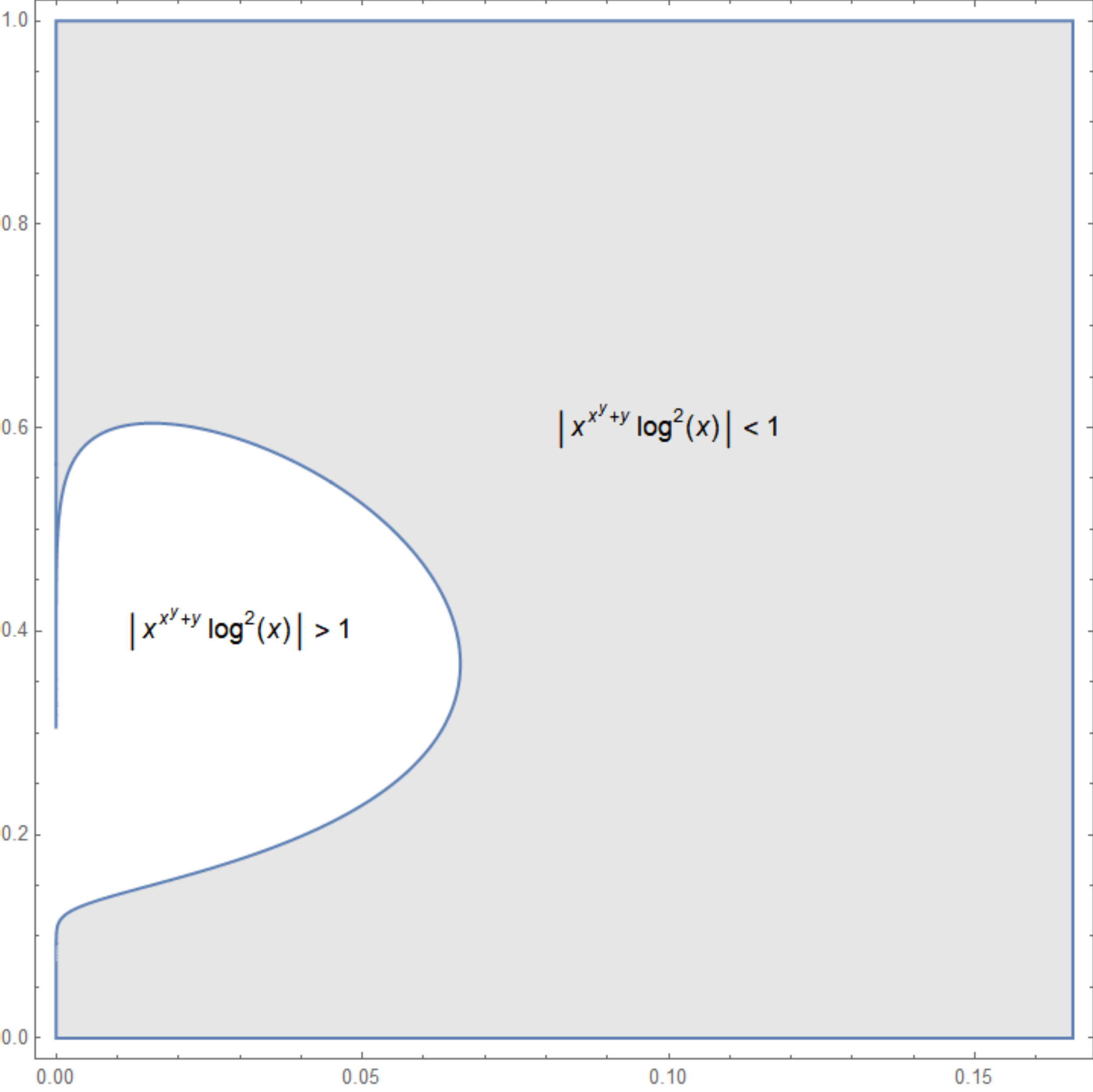}

\caption{RegionPlot of $\left| {x^{x^{y}+y}\ln^{2}x} \right|<1$}
\label{fig16}
\end{figure}

Anyway, using the \textit{RegionPlot} command of \textit{Mathematica} we can visualize it 
(\figurename~\ref{fig16}). The double step sequence converges in 
the gray region and does not in the white one.

In the gray region, where the double step sequence converge, it will 
converge to the sequence whose fixed points are given by the transcendental 
equation 
\[
y=x^{x^{y}}
\]

Now we can now put all the pieces together.

\begin{figure}[H]
\centering\includegraphics[width=360pt,scale=1]{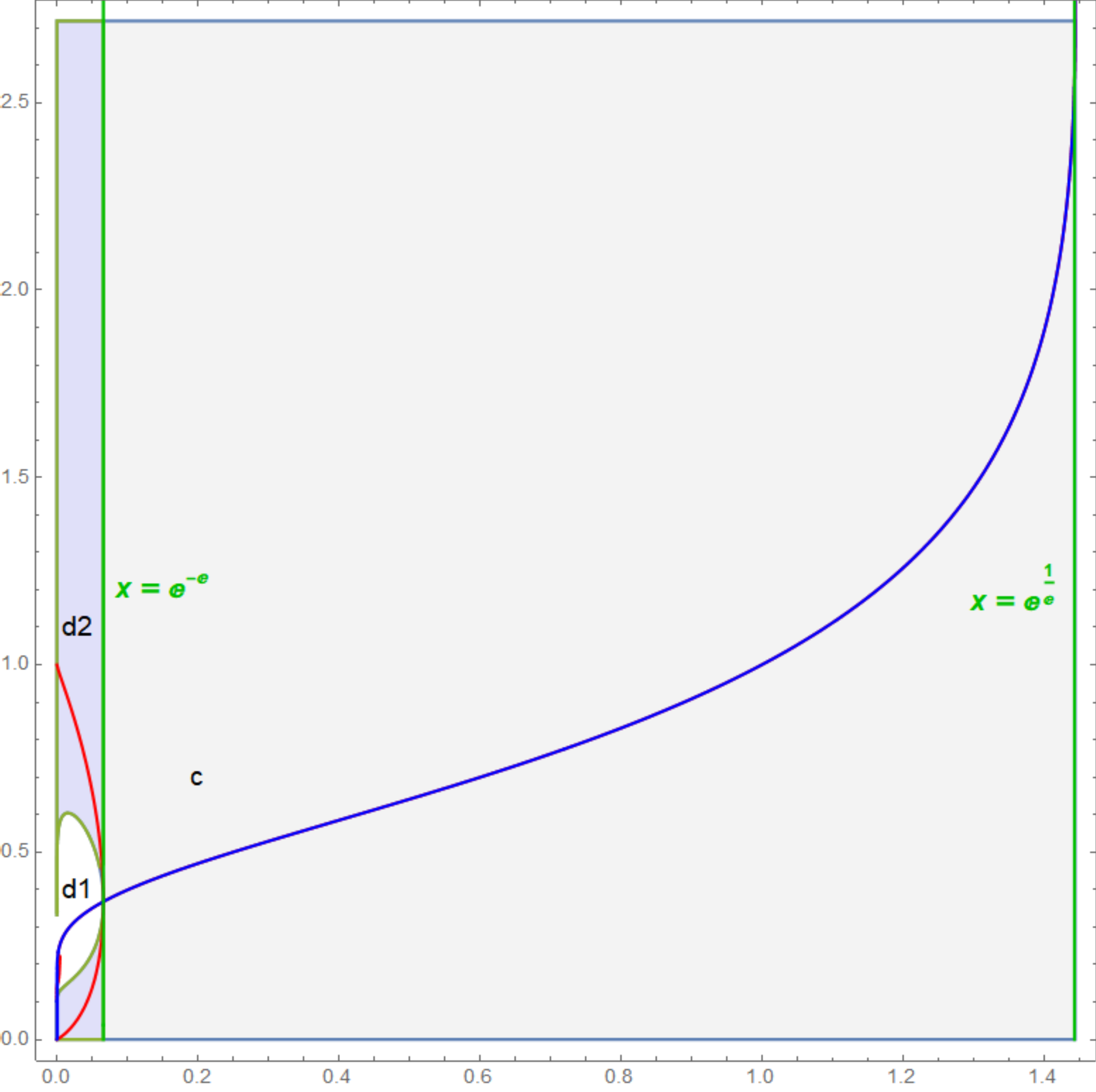}
\caption{All the plots together: $y=x^y$ and $y=x^{x^{y}}$}
\label{fig17}
\end{figure}

In \figurename~\ref{fig17} we can see the plots $p_{1}$ and 
$p_{2}$ defined by the equations $y=x^{y}$ and $y=x^{x^{y}}$respectively. 
These equations are also the equations defining the fixed points of the 
sequences $s_{1} :y_{n+1} =x^{y_{n} }$ and $s_{2} :y_{n+2} =x^{x^{y_{n} }}$. 
The gray region is where both sequences converge (``c'' in the figure), 
while the white area (``d1'' in the figure) is a region where there's no 
convergence. The blue line is produced by both equations (since the fixed 
points of a ``single iteration'' sequence are also fixed points for the one 
with a ``double iteration'' step). At the left of the line $x=e^{-e}$ 
(``d2'' in the figure) there's no convergence for $s_{1} $ and we have three 
branches. The upper and lower branches (in red) are produced only by the 
equation $y=x^{x^{y}}$.\\
Since they lie in a region of convergence for this 
sequence their values can be also produced by the infinite power tower 
function and we'll have alternating values, one on the upper branch (for 
even heights of the tower) and the other on the lower branch (odd values of 
the heights). The middle branch represent points produced by both equations. 
Anyway this branch is entirely located in the region ``d1'' where there is no 
convergence for both $s_{1} $ and $s_{2} $. The infinite power tower won't 
assume these values.

\begin{figure}[H]
\centering\includegraphics[width=360pt,scale=1]{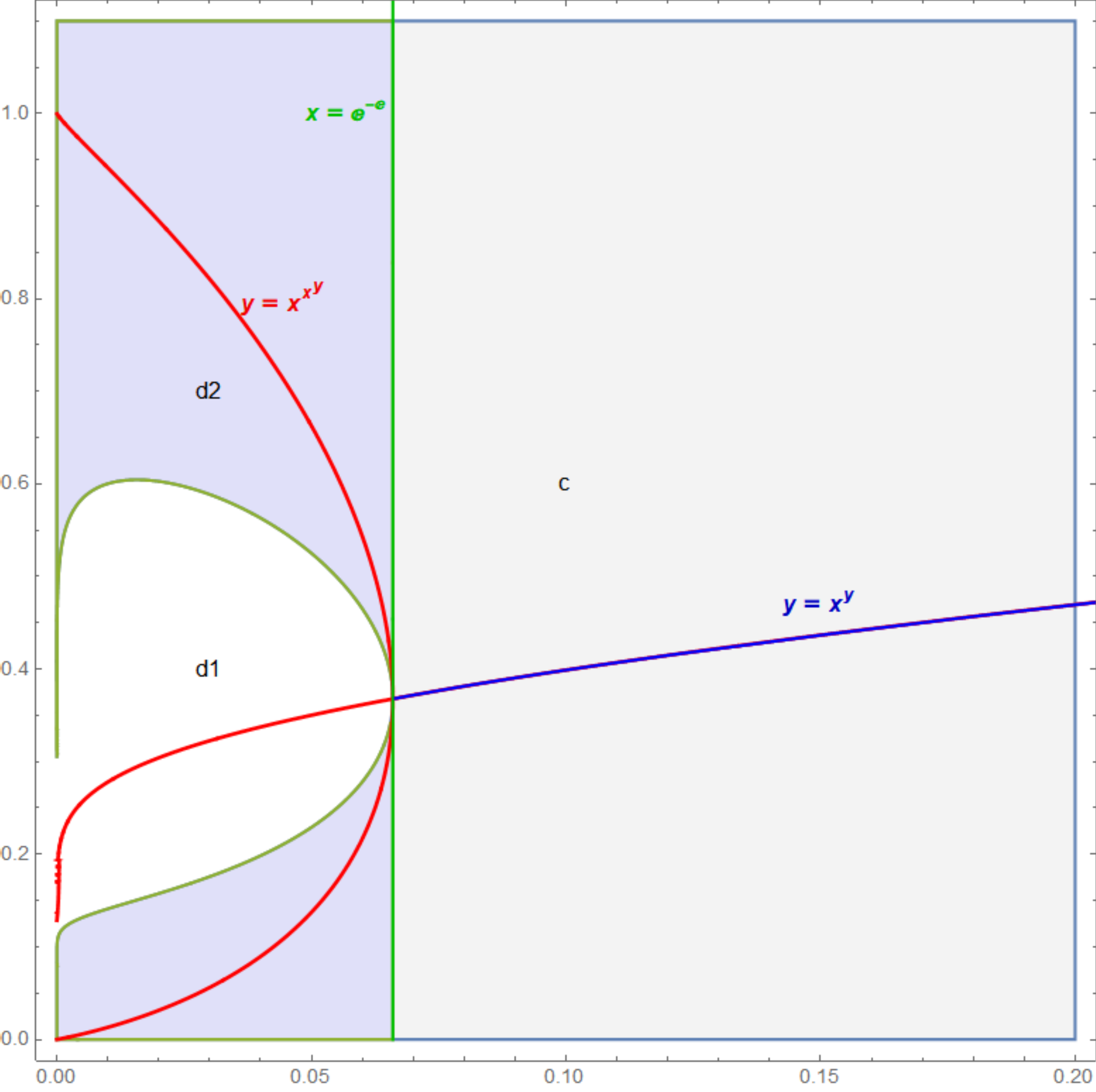}
\caption{Detail of  $y=x^y$ and $y=x^{x^{y}}$ near the “pitchfork bifurcation”}
\label{fig18}
\end{figure}

\figurename~\ref{fig18} shows an enlargement of the region with 
the three branches.

Now we can find an answer to our previously unanswered question: what is the 
limit of the infinite power tower function when $x\to 0$? 

The answer is: that limit doesn't exist. In fact, more precisely, we can have two distinct values for that limit.

That's because, since in the converging region it is $y=x^{x^{y}}$, we have 

\[
{\begin{array}{*{20}l}
 {\mbox{for\, }y\to 1} & {\lim\limits_{{\scriptsize \begin{array}{l}
 y\to 1 \\ 
 x\to 0 \\ 
 \end{array}}} x^{x^{y}}=\lim\limits_{x\to 0} x^{x}=\lim\limits_{x\to 
0} e^{x\ln x}=1} \\
 {\mbox{for\, }y\to 0} & {\lim\limits_{{\scriptsize \begin{array}{l}
 y\to 0 \\ 
 x\to 0 \\ 
 \end{array}}} x^{x^{y}}=\lim\limits_{x\to 0} x^{1}=\lim\limits_{x\to 
0} x=0} \\
\end{array} }
\]
and the equation $y=x^{x^{y}}$is verified for both $\left\{ 
{{\begin{array}{*{20}c}
 {x\to 0} \\
 {y\to 1} \\
\end{array} }} \right.$ and $\left\{ {{\begin{array}{*{20}c}
 {x\to 0} \\
 {y\to 0} \\
\end{array} }} \right.$

The emergence of the transition from a single fixed point to a 2-cycle can 
be better understood by seeing how the function $y_{n+2} =x^{x^{y_{n} }}$ 
changes with different values of the $x$. Again, let's use $z_{2} =y_{n+2} $ 
and $y=y_{n} $ 

For $e^{-e}\le x\le e^{1 \mathord{\left/ {\vphantom {1 e}} \right. 
\kern-\nulldelimiterspace} e}$ our double step function $z_{2} =x^{x^{y}}$ 
(solid line in the figures) has the same fixed points of $z=x^{y}$ (dashed 
line). Anyway, when $x<e^{-e}$, two new intersections with the identity line 
appear \figurename~(\ref{fig22}). They correspond to the values of 
the stable 2-cycle. In the meantime, the fixed point $y^{\ast }_{1} $ change 
from attractive to repulsive. In the theory of dynamical systems the 
transition from one fixed point to three fixed points is called 
\textbf{\textit{pitchfork bifurcation}}.

\begin{figure}[H]
	\hspace{-30pt}
	\begin{minipage}[b]{0.50\linewidth}
		\centering
		\includegraphics[width=0.9\linewidth, scale=1]{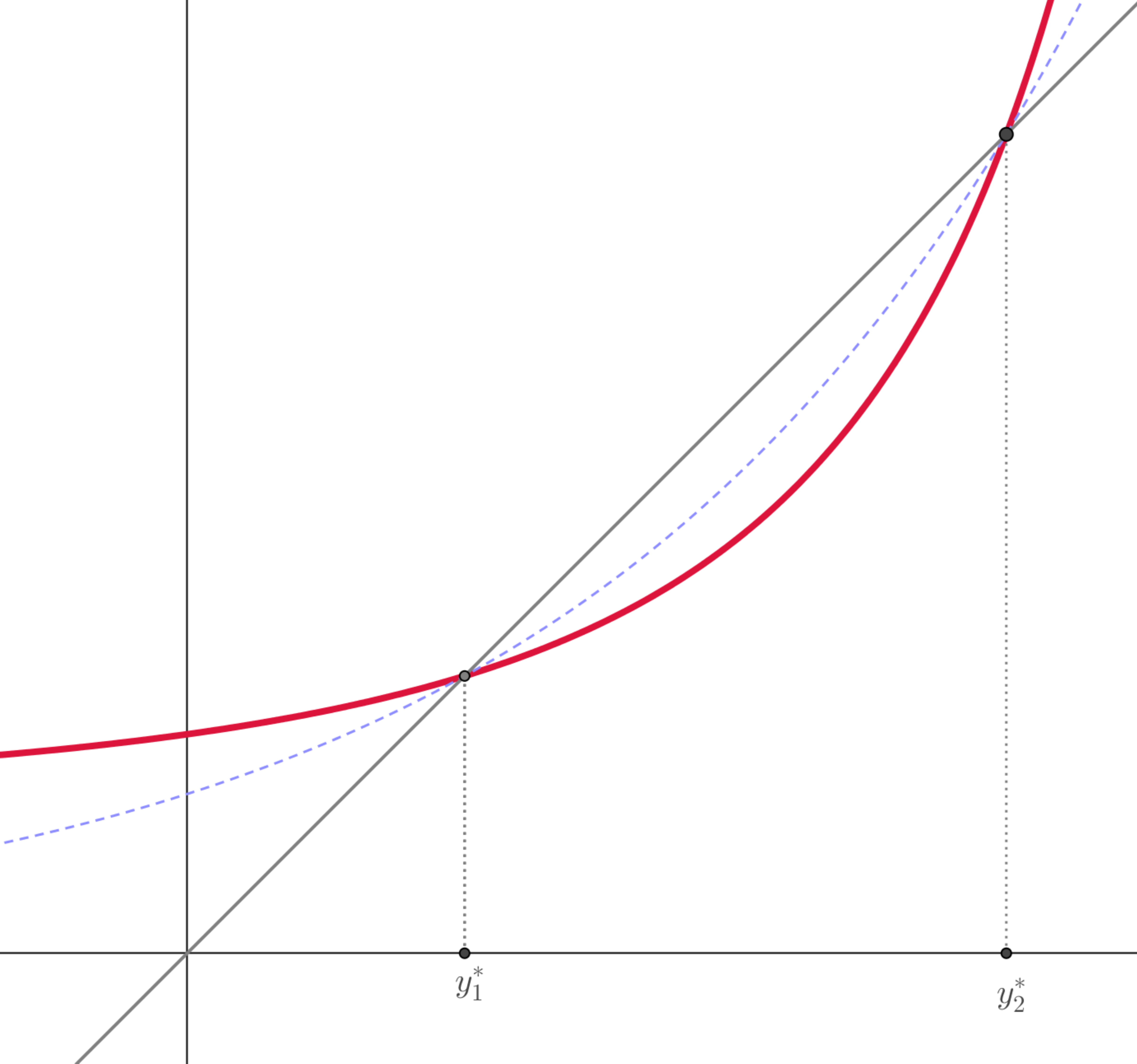}
		\caption{Graph of $z_{2}=x^{{x}^{y}}$ for $1<x<e^{1/e}$}
		\label{fig19}
	\end{minipage}\hspace*{5pt}
	\begin{minipage}[b]{0.50\linewidth}
		\centering
		\includegraphics[width=0.9\linewidth, scale=1]{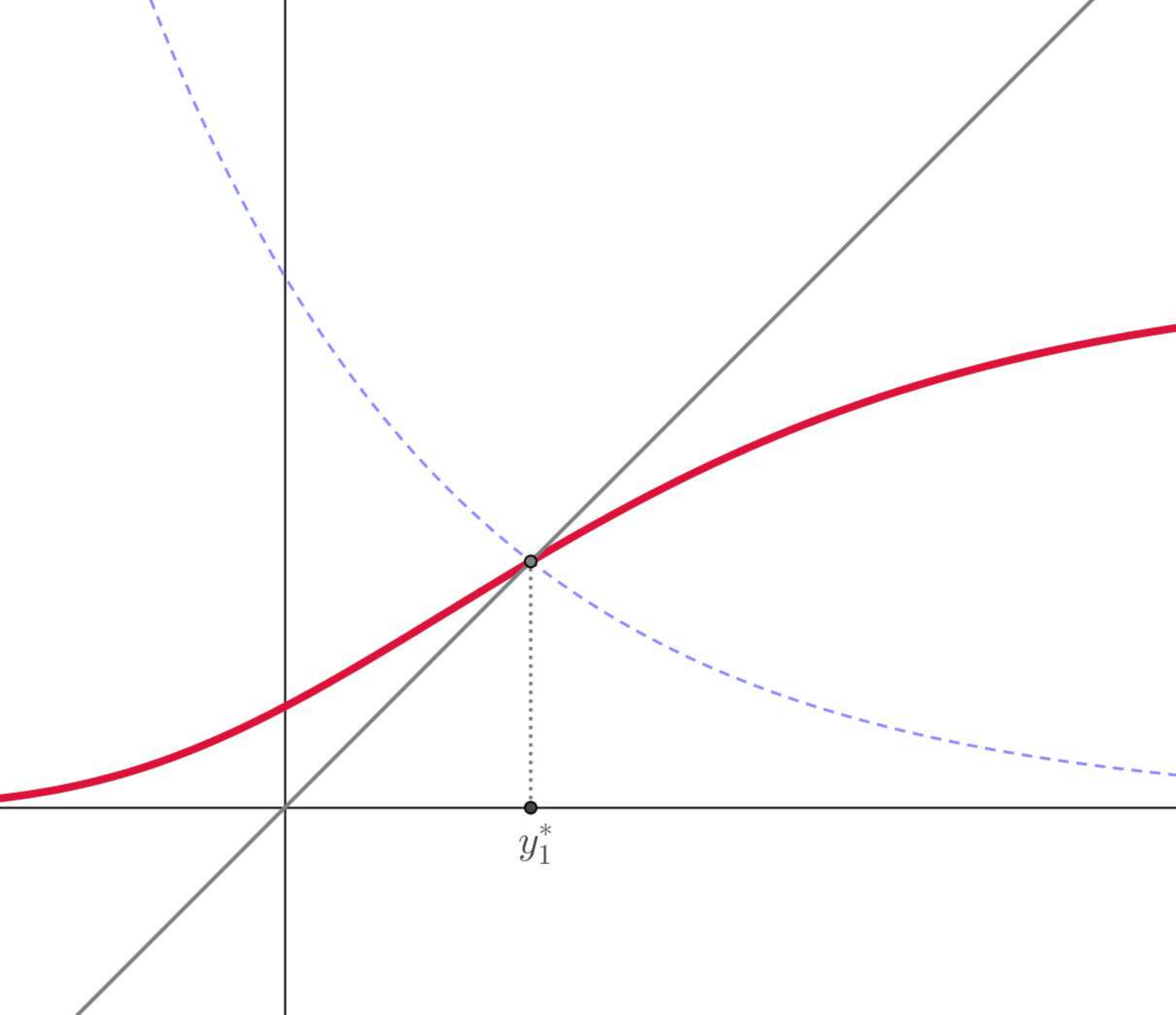}
		\caption{Graph of $z=x^{x^{y}}$ for $e^{-e}<x<1$.}
		\label{fig20}
	\end{minipage}
\end{figure}

\begin{figure}[H]
	\hspace{-30pt}
	\begin{minipage}[b]{0.5\linewidth}
		\centering
		\includegraphics[width=0.9\linewidth, scale=1]{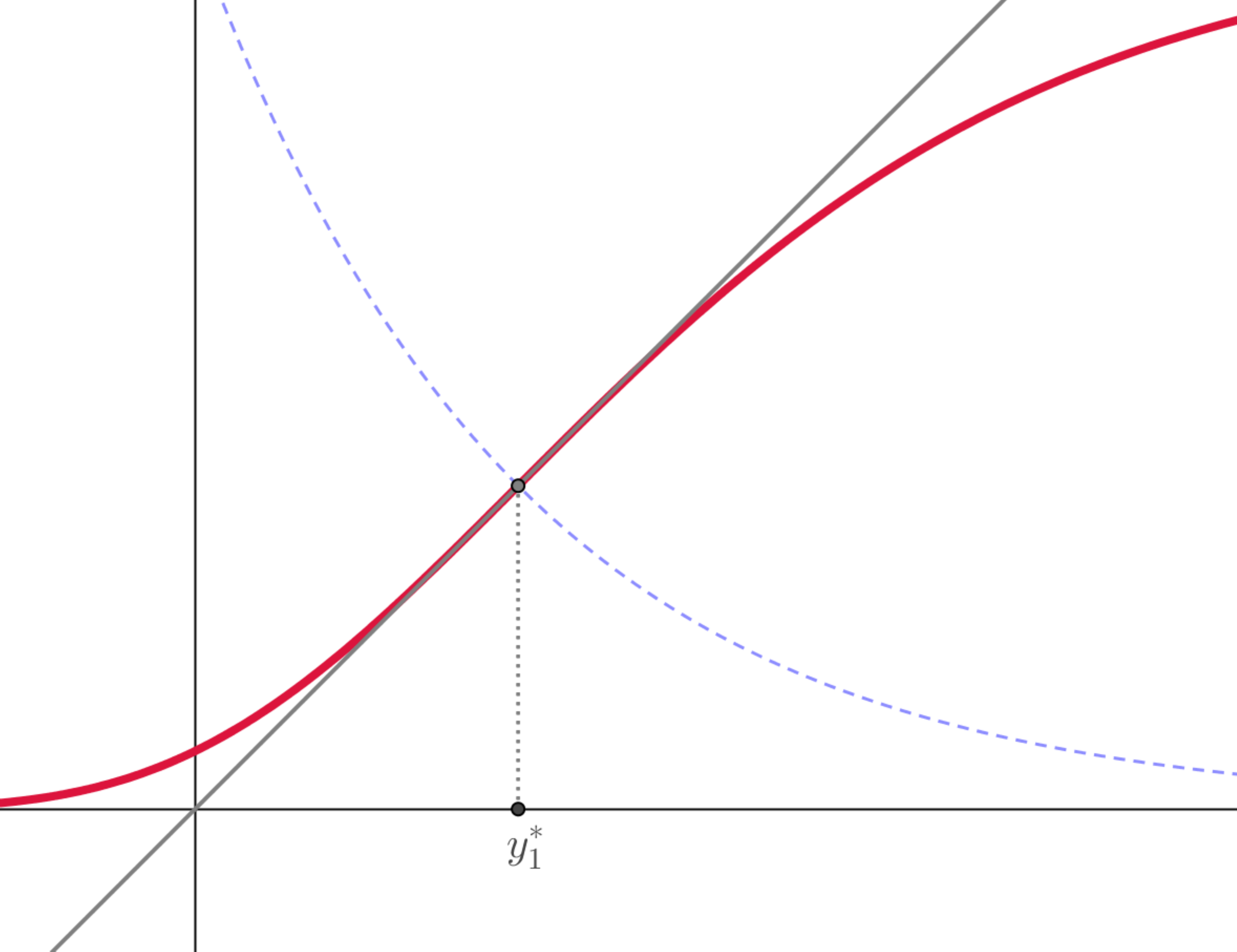}
		\caption{Graph of $z=x^{x^{y}}$ for $x=e^{-e}$.}
		\label{fig21}
	\end{minipage}\hspace*{10pt} 
	\begin{minipage}[b]{0.5\linewidth}
		\centering
		\includegraphics[width=0.9\linewidth, scale=1]{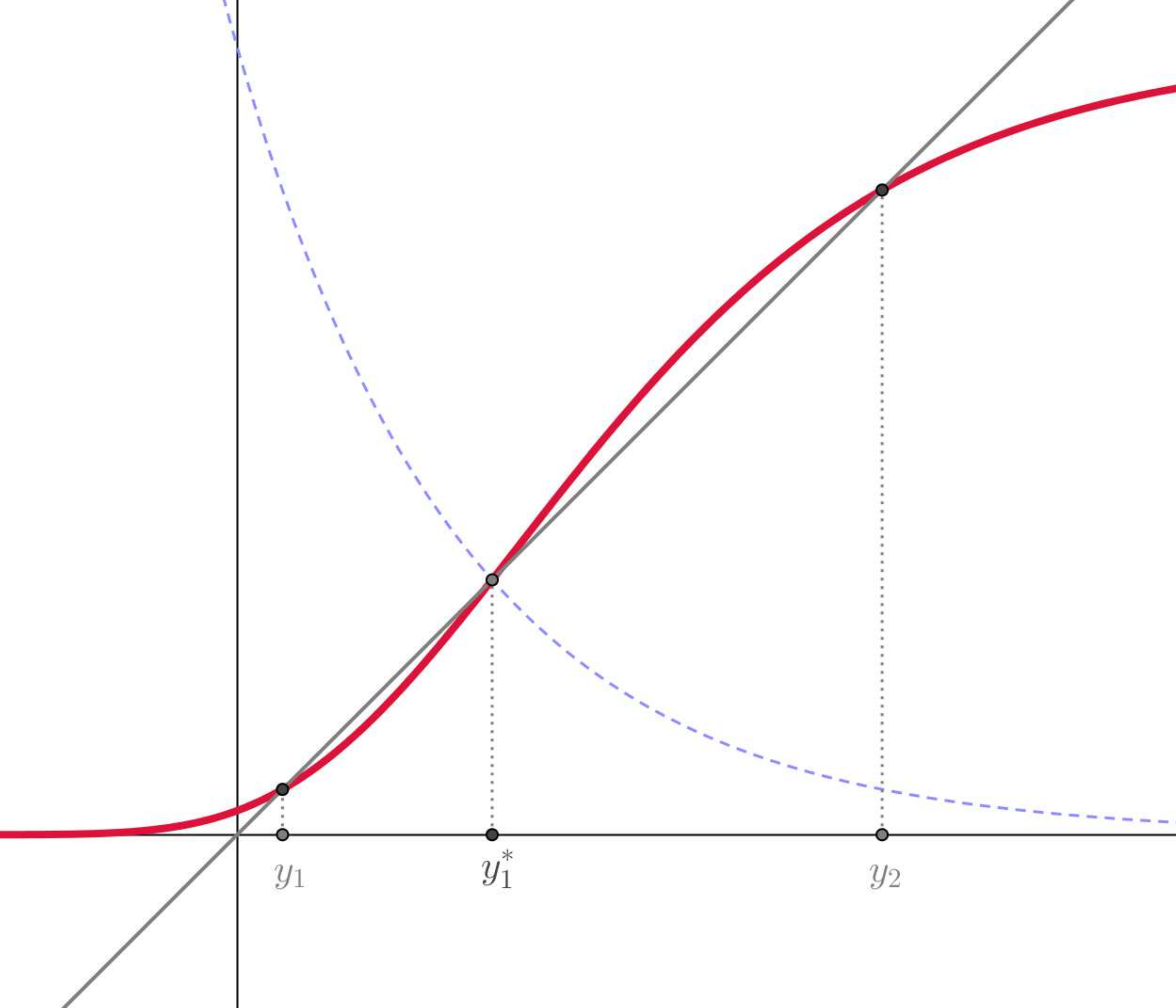}
		\caption{Graph of $z=x^{x^{y}}$ for $0<x<e^{-e}$.}
		\label{fig22}
	\end{minipage}
\end{figure}

\newpage 

\section{Some history about the power tower }
\label{sec:history}
What is the origin of the power tower function? How come that someone had 
the idea of creating such a \textit{monster}? Actually its genesis can, somehow, be 
connected with the arithmetical operations based on Peano's axioms\footnote{ 
http://mathworld.wolfram.com/PeanosAxioms.html}:
\begin{enumerate}
\footnotesize
\item zero is a number.
\item if $a$ is a number, the successor $S\left( a \right)$ of $a$ is a number.
\item zero is not the successor of a number.
\item two numbers of which the successors are equal are themselves equal.
\item If a set $K$ of numbers contains zero and also the successor of every number in $K$, then every number is in $K $(induction axiom).
\end{enumerate}
Peano's axioms are the basis of the arithmetic of natural numbers, where the 
operations of addition, multiplication and exponentiation can be defined. 
Yet the only (unary) operation included in Peano's axiom is the successor.

However, we can build the other operations by iterating the one defined at 
the previous step. The operations defined in this way are called 
\textit{hyperoperations}, and the grade 0 of this sequence is the \textit{successor} operation that, if iterated, can 
be used to define any natural number.

So we can build the sequence of operations shown in the following table: 

\begin{table}[H]
\footnotesize
\begin{center}
\begin{tabular}{|m{40pt}|m{80pt}|m{125pt}|}
\hline
& 
\vspace{4pt} \textbf{Name}& 
\vspace{4pt} \textbf{Definition} \\
\hline
\hline
\vspace{4pt} hyper0& 
\vspace{4pt} Successor & 
\vspace{4pt} $S\left( n \right)=n+1$ \\[8pt]
\hline
\vspace{4pt} hyper1& 
\vspace{4pt} Addition& 
\vspace{4pt} $n+m=S^{m}\left( n \right)$ \\[8pt]
\hline
\vspace{4pt} hyper2& 
\vspace{4pt} Multiplication& 
\vspace{4pt} \vspace{4pt} $n\cdot m=\underbrace {n+n+n+...+n}_m$ \\
\hline
\vspace{4pt} hyper3& 
\vspace{4pt} Exponentiation& 
\vspace{4pt} $n^{m}=\underbrace {n\cdot n\cdot n\cdot ...\cdot n}_m$ \\
\hline
\vspace{4pt} hyper4& 
\vspace{4pt} Tetration& 
\vspace{4pt} ${ }^{m}n=\underbrace {n^{n^{n^{...^{n}}}}}_m$ \\
\hline
\end{tabular}
\label{tab6}
\end{center}
\end{table}
\vspace{-16pt}
The sequence of hyperoperations can go on with the hyper5 (pentation), the 
hyper6 (hexation) and beyond.

Naturally, the commonly used operations are the ones reaching hyper3 
(exponentiation), but we can see that the \textit{tetration} is not just an exotic oddity but 
can be thought of as an extension of the process leading to the most usual 
arithmetical operations. 

The tetration with infinite height (infinite power tower) is often dealt 
together with the \textit{Lambert W} function (called \texttt{ProductLog} in \textit{Mathematica} and \texttt{LambertW} 
in \textit{Geogebra}).

The Lambert W function $y=W\left( x \right)$ is defined as the inverse 
function of $x=y\cdot e^{y}$ (note that there's no algebraic closed form 
expression for this function). \\
The LambertW function can be used to solve 
certain type of transcendental equations such as, for instance, 
$x\,e^{x}=2$. Its solution can be written as $x=LambertW\left( 2 \right)$ 
and the numerical value returned is $0.852606$ (since $0.852606\cdot 
e^{0.852606}=2)$.

\begin{minipage}[m]{0.75\linewidth}
Taking advantage of the definition of the LambertW function, the fixed 
points of the infinite power tower can be expressed as
\[
y=x^{y}=\frac{W\left( {-\ln x} \right)}{-\ln x}
\]
In fact, starting from

$y=x^{y}\to y=e^{y\ln x}\to ye^{-y\ln x}=1,$ multiply both sides by $-\ln x$
\[
-y\ln x\cdot e^{-y\ln x}=-\ln x
\]
set $w=-y\ln x;	z=-\ln x$

$we^{w}=z\to w=W\left( z \right)$ that is, by definition, the LambertW 
function.

Substitute back the $w$ and $z$

$-y\ln x=W\left( {-\ln x} \right)\to y=\frac{W\left( {-\ln x} \right)}{-\ln 
	x}$ that is the explicit form of the implicit function defined by $y=x^{y}$
\end{minipage}\hspace*{10pt} 
\begin{minipage}[m]{0.25\linewidth}
	\centering
	\includegraphics[width=0.9\linewidth, scale=1]{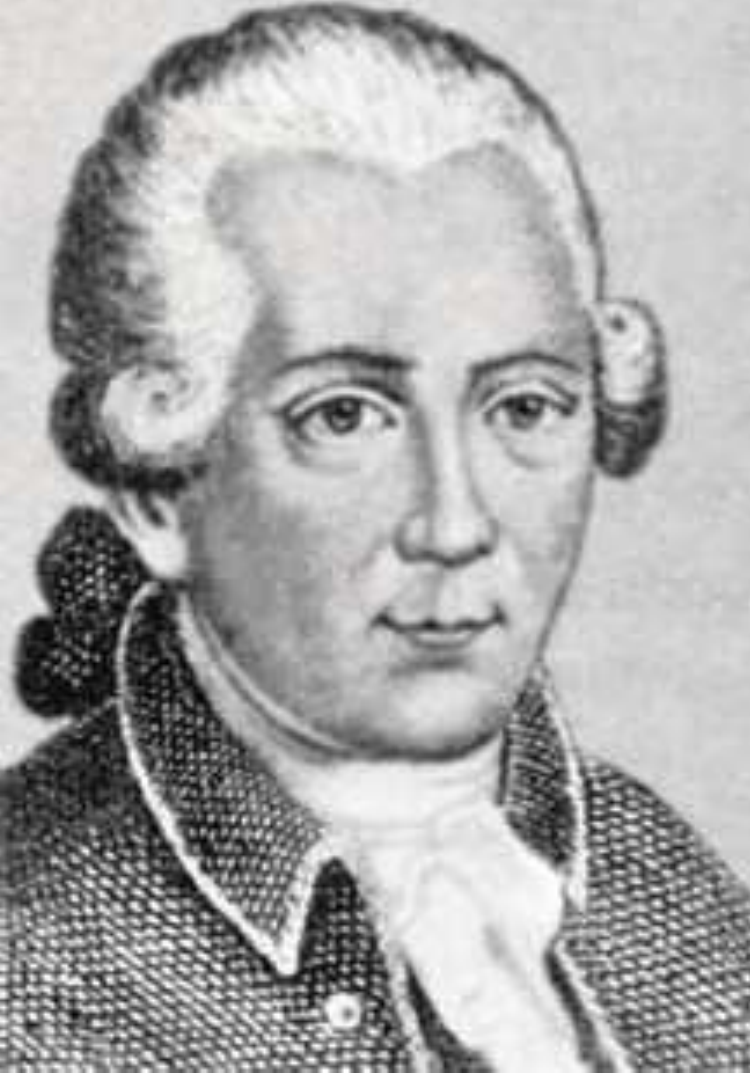}
	\footnotesize Johann Heinrich Lambert (1728-1777)
\end{minipage}

The definition of the Lambert W function originated by the article 
``\textit{Observationes variae in mathesin puram''}\footnote{ Lambert J. H., (1758). \textit{Observationes variae in mathesin puram}, Acta Helveticae 
physico-mathematico-anatomico-botanico-medica, Band III, 128--168, 
1758.\underline { }} published in 1758 by the Swiss mathematician 
\textbf{Johann Heinrich Lambert} in which he dealt with the solution of the 
trinomial transcendental equation $x^{m}+px=q$ and discovered that, under 
certain conditions, the solution (a solution) could be expressed with the 
following series: 

\scriptsize{
\[
x=\frac{q}{p}-\frac{q^{m}}{p^{m+1}}+m\frac{q^{2m-1}}{p^{2m+1}}-m\frac{3m-1}{2}\frac{q^{3m-2}}{p^{3m+1}}+m\frac{4m-1}{2}\,\,\frac{4m-2}{3}\frac{q^{4m-3}}{p^{4m+1}}-m\frac{5m-1}{2}\,\,\frac{5m-2}{3}\,\,\frac{5m-3}{4}\frac{q^{5m-4}}{p^{5m+1}}
\]
}
\normalsize
To derive above series Lambert used a procedure that was later generalized 
by \textbf{Joseph-Louis Lagrange} in 1770\footnote{ Lagrange, Joseph-Louis 
(1770). \textit{Nouvelle m\'{e}thode pour r\'{e}soudre les \'{e}quations litt\'{e}rales par le moyen des s\'{e}ries}, M\'{e}moires de l'Acad\'{e}mie Royale des Sciences et 
Belles-Lettres de Berlin. 24: 251--326.} with what's presently known as 
``\textit{Lagrange inversion theorem}''.

With Lagrange's method, given a polynomial function\footnote{ If the polynomial contains a constant 
term $a_{0} $ it's possible to eliminate it with a change of variable 
$y\mapsto y-a_{0} $ } $y=f\left( x \right)=a_{1} x+a_{2} x^{2}+...+a_{m} 
x^{m}$ it's possible to find the series expansion of the inverse function 
$x=g\left( y \right)=A_{1} y+A_{2} y^{2}+...$ by applying the following 
steps\footnote{ http://mathworld.wolfram.com/SeriesReversion.html}:

$\bullet $ plug the first expression in the second
\[
x=A_{1} \left( {a_{1} x+a_{2} x^{2}+...+a_{m} x^{m}} \right)+A_{2} \left( 
{a_{1} x+a_{2} x^{2}+...+a_{m} x^{m}} \right)^{2}+A_{3} \left( {a_{1} 
x+a_{2} x^{2}+...+a_{m} x^{m}} \right)^{3}+...
\]
$\bullet $ equate the coefficients of the right and left sides having the same grade of the $x$.

\[
{\begin{array}{*{20}c}
 {A_{1} a_{1} =1} & \to & {A_{1} =\frac{1}{a_{1} }} \\
 {A_{1} a_{2} +A_{2} a_{1}^{2}=0} & \to & {A_{2} =-\frac{a_{2} }{a_{1} 
^{3}}} \\
 {A_{1} a_{3} +2A_{2} a_{1} a_{2} +A_{3} a_{1}^{3}=0} & \to & {A_{3} 
=\frac{2a_{2}^{2}-a_{1} a_{3} }{a_{1}^{5}}} \\
 {...} & & {...} \\
\end{array} }
\]

\begin{minipage}[m]{0.75\linewidth}
By finding the inverse function (or, better, an approximation of the inverse 
function around the point $x_{0} =0)$ it is also possible to find the 
approximate value of a root $\tilde{{x}}$ of a polynomial equation having 
the form $f\left( x \right)=q$ since it is $\tilde{{x}}=g\left( {f\left( 
	{\tilde{{x}}} \right)} \right)=g\left( q \right)$.

This procedure can be extended to generic (not polynomial) functions 
$z=f\left( w \right)\to w=g\left( z \right)$ using a more general form of 
the Lagrange inversion theorem. Naturally there is the problem of 
convergence of the series, problem that we won't discuss here.
\end{minipage}\hspace*{10pt} 
\begin{minipage}[m]{0.25\linewidth}
	\centering
	\includegraphics[width=0.9\linewidth, scale=1]{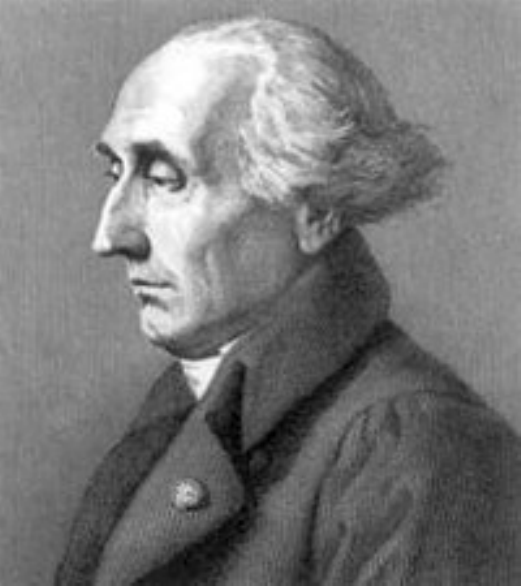}
	\footnotesize Joseph-Louis Lagrange (1736-1813)
\end{minipage}

In a subsequent article, ``\textit{Observations analytiques''}$^{\, }$\footnote{ Lambert J. H., (1770). 
\textit{Observations analytiques,} Nouveaux M\'{e}moires de l'Acad\'{e}mie royale des sciences de Berlin, 
ann\'{e}e 1770/1772} published in 1772, Lambert, examined the similar 
trinomial equation $x=q+x^{m}$ and wrote down the series that express not 
only a root of the equation, but also the powers of that root. In this 
article Lambert also mentions his meeting with L. Euler in Berlin in 1764 
and their discussions about the series connected with polynomial equations.

Some years later, in 1779, \textbf{Leonhard Euler} published ``\textit{De serie Lambertina plurimisque eius insignibus proprietaribus}''\footnote{ 
	Euler L., (1779). \textit{De serie Lambertina plurimisque eius insignibus proprietaribus, } originally published in ``Acta Academiae Scientarum 
	Imperialis Petropolitinae'' 1779, 1783, pp. 29-51} in which, referring to 
the previous works by Lambert, he investigated the solutions of another 
trinomial equation, equivalent to the one studied by Lambert, having the 
form 
\[
x^{\alpha }-x^{\beta }=v\left( {\alpha -\beta } \right)x^{\alpha +\beta }
\]
The equivalence can be verified by choosing the transformation of the 
parameters $\alpha =-m$, $\beta =-1$, $v\left( {\alpha -\beta }
\right)=q$, obtaining
\[
\frac{1}{x^{m}}-\frac{1}{x}=\frac{v\left( {\alpha -\beta } 
	\right)}{x^{m+1}}\to x-x^{m}=q
\]

In this case the series useful to express the solution (or one ot its powers) is\footnote{ Euler's 
series is not equivalent to Lambert's because Euler's series is centered 
around 1 and Lambert's is centered around 0. }: 
\scriptsize
\[
x^{n}=1+nv+\frac{1}{2}n\left( {n+\alpha +\beta } 
\right)v^{2}+\frac{1}{6}n\left( {n+\alpha +2\beta } \right)\left( {n+2\alpha 
+\beta } \right)v^{3}+\frac{1}{24}n\left( {n+\alpha +3\beta } \right)\left( 
{n+2\alpha +2\beta } \right)\left( {n+3\alpha +\beta } \right)v^{4}+...
\]
\normalsize
Euler then makes a transformation of both expressions in the special cases 
$\alpha \to \beta \to 1$ for the first equation and $\alpha \to \beta \to 
1\wedge n\to 0$ for the second.

For the first expression ($x^{\alpha }-x^{\beta }=v\left( {\alpha -\beta } 
\right)x^{\alpha +\beta })$ it is 
\[
\frac{x^{\alpha }-x^{\beta }}{\alpha -\beta }=vx^{\alpha +\beta }\to 
\frac{x^{\alpha }\left( {1-x^{\beta -\alpha }} \right)}{\alpha -\beta 
}=vx^{\alpha +\beta }\to \frac{x^{\alpha }\left( {x^{\beta -\alpha }-1} 
\right)}{\beta -\alpha }=vx^{\alpha +\beta }
\]
and taking the limit $\beta \to \alpha ,\,\,\,\,\,\beta -\alpha \to 
\varepsilon ,\,\,\,\,\,\,\varepsilon \to 0$ 
\[
\lim\limits_{\varepsilon \to 0} \frac{x^{\alpha }\left( {x^{\varepsilon 
}-1} \right)}{\varepsilon }=vx^{\alpha +\alpha }\Rightarrow x^{\alpha }\ln 
x=vx^{2\alpha }\Rightarrow \ln x=vx^{\alpha }
\]
and, since $\alpha \to 1$, $\ln x=vx$.

For the second expression it is
\scriptsize
\[
\frac{x^{n}-1}{n}=v+\frac{1}{2}\left( {n+\alpha +\beta } 
\right)v^{2}+\frac{1}{6}\left( {n+\alpha +2\beta } \right)\left( {n+2\alpha 
+\beta } \right)v^{3}+\frac{1}{24}\left( {n+\alpha +3\beta } \right)\left( 
{n+2\alpha +2\beta } \right)\left( {n+3\alpha +\beta } \right)v^{4}+...
\]
\normalsize
and taking the limits for $\alpha \to 1,\,\,\,\beta \to 
1,\,\,\,n\to 0\,\,$ it is
\[
\lim\limits_{n\to 0} \frac{x^{n}-1}{n}=\ln 
x=v+v^{2}+\frac{3}{2}v^{3}+\frac{8}{3}v^{4}+...
\]
\vspace{-20pt}
\begin{minipage}[m]{0.75\linewidth}
Putting together the two expressions we can say that a special solution of 
Euler's trinomial equation can be written in two different ways:\\
(1) the 
solution of $\ln x=vx$ and \\
(2) the result of the series $\ln 
x=v+v^{2}+\frac{3}{2}v^{3}+\frac{8}{3}v^{4}+...$

This means that the solution of the transcendental equation $\ln x=vx$ can 
be expressed by the series $\ln 
x=v+v^{2}+\frac{3}{2}v^{3}+\frac{8}{3}v^{4}+...$ and if we set $\ln x=t$ 
(and $x=e^{t})$ we have $t=ve^{t}$ whose solution is 
$t=v+v^{2}+\frac{3}{2}v^{3}+\frac{8}{3}v^{4}+...$

The equation $t=ve^{t}$ can be rewritten as $-te^{-t}=-v$ and, using the 
definition of the Lambert W function as solution of $we^{w}=z\Rightarrow 
w=W\left( z \right)$ we have $-t=W\left( {-v} \right)$ that is $t=-W\left( 
{-v} \right)$ . So here we have the series expansion 
\end{minipage}\hspace*{10pt} 
\begin{minipage}[m]{0.25\linewidth}
	\centering
	\includegraphics[width=0.9\linewidth, scale=1]{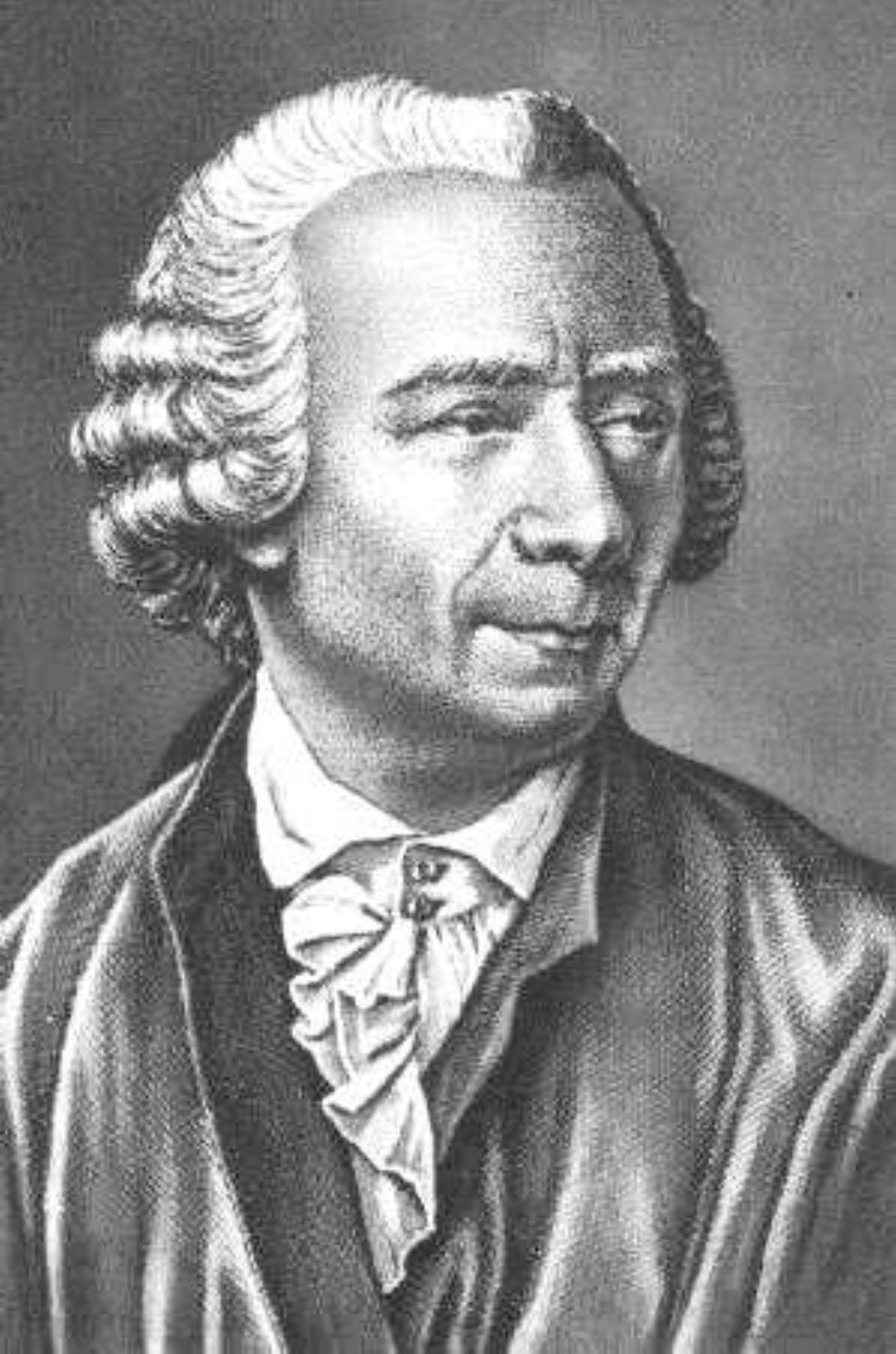}
	\footnotesize Leonhard Euler (1707-1783))
\end{minipage} 
$-W\left( {-v} \right)=v+v^{2}+\frac{3}{2}v^{3}+\frac{8}{3}v^{4}+...$ and 
$W\left( {-v} \right)=-v-v^{2}-\frac{3}{2}v^{3}-\frac{8}{3}v^{4}+...$

that is (setting $-v=z)$, \[W\left( z 
\right)=z-z^{2}+\frac{3}{2}z^{3}-\frac{8}{3}z^{4}+...\]

representing the series expansion of the LambertW function.

Even more closely related with the subject of this article is another work by Euler: ``\textit{De formulis exponentialibus replicatis''}\footnote{ Euler L., (1777). \textit{De formulis exponentialibus replicatis,} presented to the St. Petersburg 
	Academy in 1777 and published in ``Acta Academiae Scientarum Imperialis 
	Petropolitinae 1, 1778''. Also in Opera Omnia: Series 1, Volume 15, pp. 268 
	-- 297.}, presented in 1777 (two years before the publication of  ``\textit{De serie Lambertina}''), in which he investigated a problem posed by the French philosopher 
and mathematician Nicolas de Condorcet (known as Marquis de Condorcet), 
regarding the convergence of the sequence  $r,\,\,r^{\alpha 
},\,\,r^{r^{\alpha }},\,\,r^{r^{r^{\alpha }}},...$ 
%
%

The article's opening is very interesting to point out Euler's keen interest 
in such expressions. Its translation goes more or less like this:
\begin{quote}
\vspace{-6pt}\textit{"The famous Marquis de Condorcet recently shared with the academy deep 
speculations regarding some rather unfamiliar analytic formulas, among which 
we can, first of all, include the formulas called repeated exponentiations, 
where every power goes into the power exponent following it.
Yet, little has been achieved about the nature of such expressions and 
despite the force of those investigations, led with incredible sagacity, no 
clear knowledge and perception has been reached. 
Hence it will not be useless to explain here some special properties of such 
expressions."}
\end{quote}
\vspace{-6pt}
In the article Euler proves that the sequence $r,\,\,r^{\alpha 
},\,\,r^{r^{\alpha }},\,\,r^{r^{r^{\alpha }}},...$ converges if 
$e^{-e}<r<e^{1 \mathord{\left/ {\vphantom {1 e}} \right. 
\kern-\nulldelimiterspace} e}$. 

He also notes (p. 57) that the sequence $\beta =r^{\alpha },\,\,\gamma 
=r^{\beta }=r^{r^{\alpha }},\,\,\delta =r^{\gamma }=r^{r^{r^{\alpha }}},...$ 
may produce an alternate sequence of two values. In fact, choosing $r=1 
\mathord{\left/ {\vphantom {1 {16}}} \right. \kern-\nulldelimiterspace} 
{16}$ and $\alpha =1 \mathord{\left/ {\vphantom {1 2}} \right. 
\kern-\nulldelimiterspace} 2$ we have
\[
\beta =r^{\alpha }=\left( {1 \mathord{\left/ {\vphantom {1 {16}}} \right. 
\kern-\nulldelimiterspace} {16}} \right)^{1 \mathord{\left/ {\vphantom {1 
2}} \right. \kern-\nulldelimiterspace} 2}=1 \mathord{\left/ {\vphantom {1 
4}} \right. \kern-\nulldelimiterspace} 4,\,\,\gamma =r^{\beta }=\left( {1 
\mathord{\left/ {\vphantom {1 {16}}} \right. \kern-\nulldelimiterspace} 
{16}} \right)^{1 \mathord{\left/ {\vphantom {1 4}} \right. 
\kern-\nulldelimiterspace} 4}=1 \mathord{\left/ {\vphantom {1 2}} \right. 
\kern-\nulldelimiterspace} 2,\,\,\delta =r^{\gamma }=\left( {1 
\mathord{\left/ {\vphantom {1 {16}}} \right. \kern-\nulldelimiterspace} 
{16}} \right)^{1 \mathord{\left/ {\vphantom {1 2}} \right. 
\kern-\nulldelimiterspace} 2}=1 \mathord{\left/ {\vphantom {1 4}} \right. 
\kern-\nulldelimiterspace} 4
\]
and so on, with the results assuming the alternating values $1 
\mathord{\left/ {\vphantom {1 2}} \right. \kern-\nulldelimiterspace} 2$ and 
$1 \mathord{\left/ {\vphantom {1 4}} \right. \kern-\nulldelimiterspace} 4$. 
Euler shows that, in general, this happens when $r^{\Phi }=\Psi $ and 
$r^{\Psi }=\Phi $ leading to the identity $\Phi^{\Phi }=\Psi^{\Psi }$, 
since it is
\[
r^{\Phi \cdot \Psi }=\Psi^{\Psi }\to \left( {r^{\Psi }} \right)^{\Phi 
}=\Psi^{\Psi }\to \Phi^{\Phi }=\Psi^{\Psi }
\]
Now, this equation doesn't necessarily imply that $\Phi =\Psi $, because the 
function $y=x^{x}$ has a turning point for $x=1 \mathord{\left/ {\vphantom 
{1 e}} \right. \kern-\nulldelimiterspace} e$ and some $y$ can be obtained with 
two different values of the $x$.

To find the relation between the two values satisfying the equation Euler 
sets $\Psi =p\cdot \Phi $ and finds that it must be\footnote{ We have used 
Euler's method in Section \ref{sec:fixed3}.}

$\Phi =p^{p \mathord{\left/ {\vphantom {p {\left( {1-p} \right)}}} \right. 
\kern-\nulldelimiterspace} {\left( {1-p} \right)}},\,\,\,\Psi =p^{1 
\mathord{\left/ {\vphantom {1 {\left( {1-p} \right)}}} \right. 
\kern-\nulldelimiterspace} {\left( {1-p} \right)}}$ and $r=\Phi^{1 
\mathord{\left/ {\vphantom {1 \Psi }} \right. \kern-\nulldelimiterspace} 
\Psi }\left( {=\Psi^{1 \mathord{\left/ {\vphantom {1 \Phi }} \right. 
\kern-\nulldelimiterspace} \Phi }} \right)$

Finally Euler asks himself which is the condition for this two values to 
converge to a single value. This happens for $p=1$ and it is $\Phi =\lim\limits_{p\to 1} p^{p \mathord{\left/ {\vphantom {p {\left( {1-p} 
\right)}}} \right. \kern-\nulldelimiterspace} {\left( {1-p} 
\right)}}=e^{-1},\,\,\,\Psi =\lim\limits_{p\to 1} p^{1 \mathord{\left/ 
{\vphantom {1 {\left( {1-p} \right)}}} \right. \kern-\nulldelimiterspace} 
{\left( {1-p} \right)}}=e^{-1}$

The corresponding value of $r$ is 
\vspace{-10pt}\[
r=\Psi^{1 \mathord{\left/ {\vphantom {1 \Phi }} \right. 
\kern-\nulldelimiterspace} \Phi }=\left( {1 \mathord{\left/ {\vphantom {1 
e}} \right. \kern-\nulldelimiterspace} e} \right)^{e}=e^{-e}
\]
He then concludes that the relations $r^{\Phi }=\Psi $ and $r^{\Psi }=\Phi $ 
will always yield two different values if $r<e^{-e}$.

\vspace{-15pt}\section{Conclusions}
\label{sec:conclusions}
We have considered the function based on a reiterated exponentiation 
$y=x^{x^{x^{{\mathinner{\mkern2mu\raise1pt\hbox{.}\mkern2mu 
\raise4pt\hbox{.}\mkern2mu\raise7pt\hbox{.}\mkern1mu}}}}}$ and have 
investigated its properties, finding some counterintuitive fact. During our 
journey we had to cope with the unusual definition of this function, with 
its infinite sequence of exponents piling up one over the others. To proceed 
forward and make some headway we had to use different mathematical 
arguments, such as the concept of function and inverse function, limits and 
derivatives, exponentials and logarithms, sequences, fixed points of 
recursive sequences, cobweb diagrams and others. We also used experimental 
empirical tools like complex numerical computations and graphical plots provided by mathematical software packages. At 
the end we can say that much of the properties characterizing the infinite 
power tower function and its convergence (or not) to finite values have been 
explained.

Anyway, what we are left with is a vague sense of awe and amazement in 
observing the mysterious metamorphosis of this function, from one leading to 
infinite results (as it was expected in the very early stages, before 
starting a more in-depth analysis), to one producing finite values and, 
lastly, to one undergoing some serious structural change (called bifurcation 
in the field of dynamical systems) and generating stable 2-cycles.

\newpage 

\section{References}
\label{sec:references}
Corless R.M., Gonnet G. H., Hare D. E. G., Jeffrey D. J., Knuth, D. E., 
(1996). \textit{On the Lambert W function}, Advances in Computational Mathematics 5, pp. 329-359.

Knoebel R. A., (1981). \textit{Exponentials Reiterated, }The American Mathematical Monthly Vol. 88, No. 4 
(Apr., 1981), pp. 235-252

Lynch P., \quad (2017). \textit{The Fractal Boundary of the Power Tower Function, }Proceedings of Recreational Mathematics Colloquium V - G4G, 
pp. 127-138

Lynch P., (2013). \textit{The Power Tower Function, }\\
\footnotesize{\url{https://thatsmaths.files.wordpress.com/2013/01/powertowerlambert.pdf}}
\normalsize

Anderson J., (2004). \textit{Iterated exponentials, }The American Mathematical Monthly Vol. 111, No. 8 
(Oct., 2004), pp. 668-679 

Glasscock D., \textit{Exponentiales replicatas} (talk notes),\\ 
\footnotesize{\url{http://mathserver.neu.edu/\textasciitilde dgglasscock/eulerexponent.pdf}
}
\normalsize

Strogatz S., (1994). \textit{Nonlinear dynamics and chaos}, Westview Press. Chapter 10: ``One-dimensional maps''

Tetration (wikipedia): \footnotesize{\url{https://en.wikipedia.org/wiki/Tetration}}
\normalsize

Hyperoperation (wikipedia): \footnotesize{\url{https://en.wikipedia.org/wiki/Hyperoperation}}
\normalsize

Peano's axioms (mathworld): \footnotesize{\url{http://mathworld.wolfram.com/PeanosAxioms.html}}
\normalsize

Series reversion (mathworld):\footnotesize{\url{http://mathworld.wolfram.com/SeriesReversion.html}}
\normalsize

\vspace{5pt}
\textbf{\textit{Historical papers}}

Lambert J. H., (1758). \textit{Observationes variae in mathesin puram}, Acta Helveticae 
physico-mathematico-anatomico-botanico-medica, Band III, 1758, pp. 128--168

Lagrange J. L., (1770). \textit{Nouvelle m\'{e}thode pour r\'{e}soudre les \'{e}quations litt\'{e}rales par le moyen des s\'{e}ries,} M\'{e}moires de l'Acad\'{e}mie Royale des Sciences 
et Belles-Lettres de Berlin. 24, 1770, pp. 251--326

Lambert J. H., (1770). \textit{Observations analytiques, }Nouveaux M\'{e}moires de l'Acad\'{e}mie royale des 
sciences de Berlin, ann\'{e}e 1770/1772

Euler L., (1779). \textit{De serie Lambertina plurimisque eius insignibus proprietaribus}, Acta Academiae Scientarum Imperialis Petropolitinae, 
1779, 1783, pp. 29-51

Euler L., (1777). \textit{De formulis exponentialibus replicatis,} presented to the St. Petersburg Academy in 1777 and 
published in Acta Academiae Scientarum Imperialis Petropolitinae 1, 1778.

\end{document}